\documentclass[11pt]{article}
\usepackage[a4paper,margin=27mm]{geometry}
\usepackage{amsmath,amssymb,mathtools,amsthm,mathrsfs}
\usepackage[T1]{fontenc}
\usepackage{lmodern}
\newcommand{\addresscaps}{\scshape}
\usepackage{hyperref}
\hypersetup{colorlinks=true,linkcolor=blue,citecolor=blue,urlcolor=blue}
\allowdisplaybreaks
\numberwithin{equation}{section}
\newtheorem{theorem}{Theorem}[section]
\newtheorem{lemma}[theorem]{Lemma}
\newtheorem{proposition}[theorem]{Proposition}
\newtheorem{corollary}[theorem]{Corollary}
\theoremstyle{definition}
\newtheorem{definition}{Definition}[section]
\theoremstyle{remark}
\newtheorem{remark}{Remark}[section]
\title{\normalsize\MakeUppercase{Endpoint well-posedness for the Boussinesq equation}}
\author{\normalsize\MakeUppercase{Mingjuan Chen, and Minjie Shan}}
\date{}
\begin{document}
\maketitle

\begin{abstract}
We prove small-data global well-posedness for the quadratic good Boussinesq equation on the real line and the circle at the endpoint $H^{-1/2}\times H^{-3/2}$. On each domain, a small invariant spectral set supports a continuous global group, with uniqueness in a fixed local integral class. In the periodic case, both components have zero mean and the flow is Lipschitz on bounded time intervals. The a priori bounds follow from a balanced perturbation determinant for the third-order Lax operator. Averaging two real spectral branches yields coercivity and frequency-tail control. On the line, we construct the endpoint flow by commuting Hamiltonians and Green-field estimates. On the circle, discrete spectral bounds and Kishimoto's endpoint local theory give global well-posedness.
\end{abstract}

\section{Introduction}\label{sec:introduction}
In this paper, we consider the Cauchy problem for the quadratic good Boussinesq equation
\[
u_{tt}-u_{xx}+u_{xxxx}+(u^2)_{xx}=0
\]
on $\mathbb M=\mathbb R$ or $\mathbb T=\mathbb R/(2\pi\mathbb Z)$. We work with its Hamiltonian form
\begin{equation}\label{Boussinesq}
\begin{cases}
u_t=\partial_xv,\\
v_t=\partial_x(u-u_{xx}-u^2),
\end{cases}
\qquad (u,v)|_{t=0}=(u_0,v_0),
\end{equation}
where $u$ and $v$ are real-valued and $v_0$ specifies a primitive of the initial velocity. Boussinesq equations arose in the study of bidirectional shallow-water waves \cite{r9}. The positive sign of the fourth-order term in the equation above gives the linear dispersion relation
\[
\omega(\xi)=|\xi|\langle\xi\rangle
\sim\begin{cases}|\xi|,&|\xi|\ll1,\\\xi^2,&|\xi|\gg1.\end{cases}
\]
The linear evolution is therefore wave-like at low frequencies and Schr\"odinger-like at high frequencies.

We study small initial data in the endpoint space
\[
X=H^{-1/2}(\mathbb R)\times H^{-3/2}(\mathbb R),
\qquad u_t(0)=\partial_xv_0\in H^{-5/2}(\mathbb R).
\]
For the periodic problem, we write
\[
X_{\mathbb T}=H^{-1/2}(\mathbb T)\times H^{-3/2}(\mathbb T),\qquad
X_{\mathbb T,0}=\{(u,v)\in X_{\mathbb T}:\widehat u(0)=\widehat v(0)=0\}.
\]
On the circle, both means are conserved by \eqref{Boussinesq}. The condition on $v$ fixes the periodic primitive; the condition on $u$ restricts the problem to the zero-mean flow. On the line, requiring $v_0\in H^{-3/2}$ imposes a low-frequency condition that is not implied by $u_t(0)\in H^{-5/2}$ alone. We write $X^s=H^s\times H^{s-1}$ on the spatial domain under consideration.

\subsection{Earlier results}\label{subsec:earlier-results}
The Cauchy problem for the good Boussinesq equation has been studied by energy and dispersive methods. Bona and Sachs \cite{r12} established local well-posedness for a generalized Boussinesq equation and proved nonlinear stability of solitary waves in a range of wave speeds, obtaining global solutions for data close to those waves. Linares \cite{r3} proved global existence for small initial data of a generalized equation. The sign of the linear part does not prevent nonlinear blow-up: Sachs \cite{r14} constructed solutions of the good Boussinesq equation that blow up in finite time, and Liu \cite{r15} proved instability and blow-up results for a generalized Boussinesq equation.

In the energy space, Mu\~noz, Poblete, and Pozo \cite{r13} showed that small solutions decay strongly in local energy on regions of the form $|x|<Ct/\log^2t$, without a parity assumption. Chen, Guo, and Shao \cite{r4} obtained local results in $L^2$ and $H^1$ and studied global existence, blow-up, and scattering for generalized Boussinesq equations, including large radial-data scattering in higher-dimensional defocusing cases.

Fourier restriction methods lowered the regularity required for local solutions. Fang and Grillakis \cite{r16} obtained a periodic local theory at the $L^2$ level. Farah \cite{r2} proved local well-posedness on the line for $s>-1/4$, and Farah and Scialom \cite{r17} obtained the same threshold on the circle. These two results use initial displacement in $H^s$ and initial velocity $\partial_x\psi$, with $\psi\in H^{s-1}$. Kishimoto and Tsugawa \cite{r18} improved the real-line threshold to $s>-1/2$ by reducing the equation to a quadratic Schr\"odinger equation and modifying the Fourier restriction spaces. On the circle, Oh and Stefanov \cite{r22} used a normal-form transformation to reach $s>-3/8$. Kishimoto \cite{r1} subsequently proved local well-posedness at $s=-1/2$ on both the line and the circle, with uniqueness in a modified restriction space and locally Lipschitz dependence on the initial data. He also proved ill-posedness below this threshold.

The estimates in this local theory are related to those for quadratic Schr\"odinger equations. Bourgain \cite{r21} introduced Fourier restriction spaces and used periodic Strichartz estimates to study nonlinear Schr\"odinger equations. Kenig, Ponce, and Vega \cite{r19} proved bilinear estimates for one-dimensional quadratic Schr\"odinger nonlinearities. For the equation with nonlinearity $u^2$, Bejenaru and Tao \cite{r20} constructed modified spaces to prove local well-posedness in $H^{-1}(\mathbb R)$ and ill-posedness below that regularity. We use the endpoint resolution space, the bilinear estimates, and the Duhamel extension from Kishimoto \cite{r1}. For the two-field system, the recovery of $v$ at low frequencies must still be addressed on the line.

At the endpoint considered here, the classical Hamiltonian is not defined for general data. Low-regularity conservation laws for integrable equations provide another source of a priori bounds. Killip, Vi\c{s}an, and Zhang \cite{r23} used perturbation determinants to construct conserved quantities controlling negative Sobolev norms for KdV, NLS, and mKdV on the line and the circle. Koch and Tataru \cite{r24} constructed conserved energies equivalent to $H^s$ norms for cubic NLS and mKdV when $s>-1/2$, and obtained corresponding energies for KdV when $s\ge-1$. For KdV at $H^{-1}$, Buckmaster and Koch \cite{r25} used the Miura transformation to construct weak solutions on the line and prove a priori bounds. Kappeler and Topalov \cite{r26} established global well-posedness on the circle through inverse spectral coordinates.

The real-line construction in this paper follows the commuting-flow method of Killip and Vi\c{s}an \cite{r5}, who proved global well-posedness for KdV in $H^{-1}$. They constructed regularized Hamiltonians from the spectral determinant and compared their flows with the KdV flow in Green-function coordinates. This method was used by Bringmann, Killip, and Vi\c{s}an \cite{r27} for fifth-order KdV in $H^{-1}(\mathbb R)$ and by Harrop-Griffiths, Killip, and Vi\c{s}an \cite{r28} for cubic NLS and mKdV in $H^s(\mathbb R)$, $s>-1/2$. For derivative NLS, Killip, Ntekoume, and Vi\c{s}an \cite{r30} proved equicontinuity and global well-posedness in $H^{1/6}$ under the mass restriction $\|q\|_2^2<4\pi$. Harrop-Griffiths, Killip, and Vi\c{s}an \cite{r31} removed the mass restriction from the equicontinuity result on the line. Harrop-Griffiths, Killip, Ntekoume, and Vi\c{s}an \cite{r29} then obtained global well-posedness for derivative NLS in $L^2(\mathbb R)$.

For Boussinesq equations, the associated scalar spectral operator has order three. Hirota \cite{r11} obtained explicit multisoliton solutions, and Kaup \cite{r41} developed inverse scattering for cubic spectral problems. McKean \cite{r42} studied Boussinesq dynamics on the circle using the associated spectral problem. Deift, Tomei, and Trubowitz \cite{r10} developed an inverse-scattering theory for a Boussinesq equation with the opposite dispersive sign. Charlier and Lenells \cite{r7} formulated a $3\times3$ Riemann--Hilbert problem for the good Boussinesq equation; Charlier, Lenells, and Wang \cite{r35} used that formulation to determine long-time asymptotics. Charlier and Lenells also constructed a Miura transformation from a quadratic derivative NLS equation to the good Boussinesq equation \cite{r8}, and developed inverse scattering for the water-wave Boussinesq equation with the bad sign \cite{r43}. General matrix inverse-scattering problems and boundary-value problems with $3\times3$ Lax pairs were treated by Caudrey \cite{r32} and Lenells \cite{r33}, respectively. Badanin and Korotyaev \cite{r34} studied third-order operators with three-point boundary conditions associated with Boussinesq's equation.

Perturbation determinants also apply to matrix-valued Lax systems. Klaus \cite{r6} proved $H^s$ a priori estimates, $-1/2<s<0$, for a quadratic derivative NLS equation with a $3\times3$ Lax pair, and identified its transmission coefficient with a renormalized Fredholm determinant. Gesztesy, Latushkin, and Makarov \cite{r36} related Evans functions defined by generalized Jost solutions to modified Fredholm determinants for first-order systems. Gesztesy and Makarov \cite{r37} computed determinants of matrix-valued semiseparable kernels through Volterra equations, with applications to Jost functions and inverse transmission coefficients. To construct the Boussinesq flow in $X$, we need determinant variations and Green identities for distributional coefficients, with estimates uniform in the spectral parameter.

Independently, Harrop-Griffiths \cite{r44} proved small-data global well-posedness for the good Boussinesq equation on the line in $H^s\times H^{s-1}$ for $s\ge-1/2$. His proof combines a two-field local theory with spectral a priori bounds and also yields an endpoint long-wave approximation by two counterpropagating KdV solutions. His result includes the real-line endpoint studied here after a fixed change of normalization. Our real-line proof uses commuting Hamiltonian flows, and we also treat the zero-mean periodic problem. The comparison of the well-posedness statements is given in Remark~\ref{rmk:harrop-griffiths}.

\subsection{Main results}\label{subsec:main}
We prove small-data global well-posedness in $X$ on the line and in $X_{\mathbb T,0}$ on the circle. To state the invariant domains of the flows, we introduce the spectral functional used in the a priori estimates.

Let $U=(u,v)$ and
\[
\begin{gathered}
L_0=\partial_x^3-\tfrac14\partial_x,\qquad
z_k=k^3-\tfrac14k,\quad k\ge1,\\
V_U f=\tfrac14\bigl[u\,\partial_x f+\partial_x(uf)\bigr]
+\tfrac1{4\sqrt3}vf.
\end{gathered}
\]
With $T_z=(L_0-z)^{-1/2}$, set
\[
B_z(U)=T_zV_UT_z,\qquad A_z(U)=\log\det{}_2(I+B_z(U)),
\]
\[
\begin{gathered}
s_k(U)=-\tfrac12\bigl[A_{z_k}(U)+A_{-z_k}(U)\bigr],\\
\mathscr S_R(U)=\int_R^\infty k s_k(U)\,dk,\qquad \mathscr S=\mathscr S_1.
\end{gathered}
\]
The logarithm is normalized by $A_z(0)=0$. The two spectral signs cancel the mixed term in the quadratic expansion. On a sufficiently small real ball in $X$, we prove
\[
\gamma\|U\|_X^2\le\mathscr S(U)\le\Gamma\|U\|_X^2,
\qquad 0<\gamma\le\Gamma.
\]
For sufficiently regular data, $\Phi_E^tU_0$ denotes the high-regularity solution of \eqref{Boussinesq} at time $t$, on its maximal interval of existence. Thus
\[
\frac{d}{dt}\Phi_E^tU_0
=\bigl(\partial_x v,\partial_x (u-\partial_x^2u-u^2)\bigr)\big|_{(u,v)=\Phi_E^tU_0},
\qquad \Phi_E^0U_0=U_0.
\]
Here $E$ is the Hamiltonian in Section~\ref{subsec:notation}. The notation $\Phi_E$ will also be used for the periodic high-regularity flow. The a priori estimates below imply global existence for the approximating trajectories in our theorems.

The class $\mathcal K_{\lambda_*}$ is defined using Kishimoto's local integral equation at one fixed scale $\lambda_*$; see Section~\ref{subsec:local}.

Our first result concerns the line.
\begin{theorem}\label{thm:line}
There exist $\varepsilon>0$ and a fixed local-coordinate scale $\lambda_*\ge1$ such that the open set
\begin{equation}\label{eq:legacy-1}
\mathcal D=\{U\in X\text{ real}:\|U\|_X<\varepsilon,
\quad\mathscr S(U)<\gamma\varepsilon^2/4\}
\end{equation}
supports a global group $\Phi$ with the following properties.
\begin{enumerate}
\item For every $T>0$,
\[
\Phi:\mathbb R\times\mathcal D\longrightarrow\mathcal D
\text{ is continuous},\qquad
U_0\longmapsto\Phi^\cdot U_0\in C([-T,T];X)
\text{ is continuous},
\]
\[
\Phi^{t+s}U_0=\Phi^t(\Phi^sU_0),\quad \Phi^0U_0=U_0,
\qquad \sup_{t\in\mathbb R}\|\Phi^tU_0\|_X<\varepsilon/2.
\]
\item For $k,R\ge1$ and $t\in\mathbb R$,
\[
A_{\pm z_k}(\Phi^tU_0)=A_{\pm z_k}(U_0),\qquad
\mathscr S_R(\Phi^tU_0)=\mathscr S_R(U_0).
\]
If $U_{0,n}\in\mathcal D\cap X^{11/2}$ and $U_{0,n}\to U_0$ in $X$, then
\[
\Phi_E^\cdot U_{0,n}\longrightarrow\Phi^\cdot U_0
\quad\text{in }C([-T,T];X).
\]
\item The trajectory $U=(u,v)=\Phi^\cdot U_0$ satisfies
\[
u_t=\partial_x v,\qquad v_t=\partial_x (u-\partial_x^2u)-\mathcal N(U),
\qquad \mathcal N(U)=\partial_x (u^2)_{\rm can},
\]
in distributions. Here $\mathcal N$ is the differentiated source of Lemma~\ref{lem:source}. The trajectory is unique among real $C_tX$ integral solutions in $\mathcal K_{\lambda_*}$ with initial datum $U_0$.
\end{enumerate}
\end{theorem}

\begin{remark}\label{rmk:harrop-griffiths}
Independently, Harrop-Griffiths \cite[Theorem 1.2]{r44} proved small-data global well-posedness on $\mathbb R$ in $H^s\times H^{s-1}$ for every $s\ge-\tfrac12$, with smallness measured in $H^{-1/2}\times H^{-3/2}$ and locally Lipschitz dependence on bounded time intervals. After a fixed normalization of the coefficients, this includes the real-line endpoint considered here. His proof combines a two-field local theory with spectral a priori estimates. Theorem~\ref{thm:line} gives a commuting-flow construction at that endpoint.
\end{remark}

\begin{remark}\label{rmk:line}
Data in $\mathcal D$ need not have finite classical energy. The source $\mathcal N(U)$ is defined by convergence of the differentiated nonlinearities of smooth solutions; the subscript $\mathrm{can}$ marks this canonical differentiated limit. We do not assume that the ordinary square of an arbitrary $H^{-1/2}$ distribution exists. Uniqueness is asserted in $\mathcal K_{\lambda_*}$, not among all $C_tX$ distributional solutions.
\end{remark}

The upper bound for $\mathscr S$ gives the following small-norm consequence.
\begin{corollary}\label{cor:small-data}
Choose
\[
0<r_0<\min\{\varepsilon,\tfrac{\varepsilon}{2}\sqrt{\gamma/\Gamma}\}.
\]
Then
\[
U_0\in X\text{ real},\quad\|U_0\|_X<r_0
\quad\Longrightarrow\quad U_0\in\mathcal D.
\]
In particular, every such datum has the global trajectory, approximation property, and fixed-class uniqueness of Theorem~\ref{thm:line}.
\end{corollary}
\begin{proof}
The upper spectral bound and the choice of $r_0$ give the two inequalities defining $\mathcal D$:
\[
\|U_0\|_X<\varepsilon,\qquad
\mathscr S(U_0)\le\Gamma\|U_0\|_X^2
<\Gamma r_0^2<\gamma\varepsilon^2/4.
\]
Hence $U_0\in\mathcal D$, and the conclusions follow from Theorem~\ref{thm:line}.
This proves the corollary.
\end{proof}
\begin{remark}\label{rmk:small-data}
The spectral set $\mathcal D$ is invariant. The norm ball in Corollary~\ref{cor:small-data} is contained in $\mathcal D$ but need not itself be invariant.
\end{remark}

For zero-mean periodic data, the local complex variable determines both fields by a bounded map. The discrete spectral estimate and the endpoint local theory give the following result.
\begin{theorem}\label{thm:circle}
There exist $r,c_2>0$ and one fixed local-coordinate scale $\lambda_{\mathbb T}\ge1$ such that
\begin{equation}\label{eq:periodic-main-domain}
\mathcal D_{\mathbb T}
=\{U\in X_{\mathbb T,0}:\|U\|_{X_{\mathbb T}}<r,
\quad\mathscr S^{\mathbb T}(U)<c_2r^2/8\}
\end{equation}
is open in $X_{\mathbb T,0}$, contains a neighborhood of zero, and supports a global group $\Phi_{\mathbb T}$. Here the periodic spectral functional is defined by
\begin{equation}\label{eq:periodic-spectral-definition}
\begin{gathered}
\mathscr S^{\mathbb T}_R(U)
=-\frac12\int_R^\infty k\bigl[A^{\mathbb T}_+(k;U)+A^{\mathbb T}_-(k;U)\bigr]\,dk,
\qquad \mathscr S^{\mathbb T}=\mathscr S^{\mathbb T}_1,\\
A^{\mathbb T}_\pm(k;U)=\log\det{}_2\bigl(I+B^{\mathbb T}_\pm(k;U)\bigr),
\qquad B^{\mathbb T}_\pm(k;U)=T^{\mathbb T}_{\pm z_k}V_UT^{\mathbb T}_{\pm z_k}.
\end{gathered}
\end{equation}
The multipliers $T^{\mathbb T}_z$ are the periodic realizations of $T_z$ on the full space $L^2(\mathbb T,dx/(2\pi))$, including the constant mode. The logarithmic branches vanish at $U=0$.
\begin{enumerate}
\item For every $T>0$ there is $C_T<\infty$ such that
\[
\sup_{|t|\le T}\|\Phi_{\mathbb T}^tU-\Phi_{\mathbb T}^tV\|_{X_{\mathbb T}}
\le C_T\|U-V\|_{X_{\mathbb T}},\qquad U,V\in\mathcal D_{\mathbb T},
\]
\[
\Phi_{\mathbb T}^{t+s}=\Phi_{\mathbb T}^t\Phi_{\mathbb T}^s,
\qquad \Phi_{\mathbb T}^0=\mathrm{id},\qquad
\sup_{t\in\mathbb R}\|\Phi_{\mathbb T}^tU\|_{X_{\mathbb T}}<r/2.
\]
The group preserves $A^\mathbb T_\pm(k;\cdot)$ and $\mathscr S_R^\mathbb T$ for every $k,R\ge1$.
\item For every smooth sequence $U_{0,n}\in\mathcal D_{\mathbb T}$ converging to $U_0$ in $X_{\mathbb T}$,
\[
\Phi_E^\cdot U_{0,n}\longrightarrow\Phi_{\mathbb T}^\cdot U_0
\quad\text{in }C([-T,T];X_{\mathbb T}).
\]
The limit satisfies the two-field equation with the source of Lemma~\ref{lem:periodic-source}, and is unique in $\mathcal K^{\mathbb T}_{\lambda_{\mathbb T}}$.
\item For every integer $j\ge0$,
\[
U_0\in\mathcal D_{\mathbb T}\cap X^{j-1/2}_{\mathbb T}
\quad\Longrightarrow\quad
\Phi_{\mathbb T}^\cdot U_0\in C(\mathbb R;X^{j-1/2}_{\mathbb T}).
\]
\end{enumerate}
\end{theorem}
\begin{remark}\label{rmk:circle}
The scalar initial velocity is $u_t(0)=\partial_x v_0$ and has zero mean. For each fixed $T$, the Lipschitz constant is uniform on $\mathcal D_{\mathbb T}$; it may grow with $T$. Uniqueness is in the fixed local class $\mathcal K^{\mathbb T}_{\lambda_{\mathbb T}}$ of Section~\ref{subsec:sources-short-windows-and-the-uniqueness-class}.
\end{remark}

\begin{remark}\label{rmk:local-scales}
For completeness, the local coordinates used in the uniqueness classes are
\[
w=u+i(1-\partial_x^2)^{-1}\partial_x v,\qquad
b_{\lambda,t_0}(s,y)=\lambda^{-2}w(t_0+\lambda^{-2}s,y/\lambda).
\]
Membership in $\mathcal K_{\lambda_*}$ or $\mathcal K^{\mathbb T}_{\lambda_{\mathbb T}}$ requires $b_{\lambda,t_0}$ to belong locally to $\mathcal W_p$ on the line or $\mathcal W_{p,\lambda}$ on the circle, respectively, and to satisfy the local integral equation, for every $t_0$. These restriction spaces are specified in Definition~\ref{def:Wp}. On the circle, $y\in\mathbb T_\lambda=\mathbb R/(2\pi\lambda\mathbb Z)$; the physical equation remains on the fixed circle $\mathbb T$.

Here $\lambda=\lambda_*$ on the line and $\lambda=\lambda_{\mathbb T}$ on the circle. The rescaled equation is
\[
ib_s+\partial_y^2b=\tfrac12\lambda^{-2}(b-\bar b)
-\Omega_\lambda(\Re b)^2,\qquad
\Omega_\lambda(\xi)=\frac{\lambda^2\xi^2}{1+\lambda^2\xi^2}.
\]
We first choose $\lambda$ so that the linear term with coefficient $\lambda^{-2}$ is small, and then choose the endpoint radius. This scale is fixed at every restart and does not depend on higher Sobolev norms. It is unrelated to the spectral parameter $k$.
\end{remark}

\subsection{Proof strategy}\label{subsec:strategy}
The a priori estimates use the perturbation determinant of the third-order operator $L=L_0+V_U$. For an analytic functional with no constant or linear term, write
\begin{equation}\label{eq:quadratic-notation}
F^{[2]}(U)=\tfrac12\left.\frac{d^2}{d\theta^2}F(\theta U)\right|_{\theta=0}.
\end{equation}
The quadratic part of a single branch $A_{z_k}$ is indefinite because of a mixed $u,v$ term. This term changes sign when $z_k$ is replaced by $-z_k$. Averaging the branches and integrating in $k$ gives $\mathscr S^{[2]}(U)\sim\|U\|_X^2$. The higher-order terms are bounded by $C\|U\|_X^3$ and can be absorbed for small data. The same calculation for $\mathscr S_R$ controls the Fourier tail. Conservation therefore supplies a uniform endpoint bound and a common frequency tail for the trajectories used in the approximation.

On the line, we apply the commuting-flow method of \cite{r5} with $H_k=C_ks_k$. At each fixed $k$, its Hamiltonian field $Y_k$ is analytic on $X$, and the spectral bound gives a global auxiliary flow $\Psi_k$. The difficulty is to compare these flows with the original equation uniformly as $k\to\infty$. A bound on a fixed high-regularity ball is insufficient, since the high norms of the trajectories may depend on $k$. We use the paired Green coordinates $\mathcal G=(r_\eta,s_\eta)$ and derive their evolution from a two-parameter resolvent identity. In these coordinates, we prove a weak Lipschitz estimate uniform in $k$ and a strong approximation estimate on each fixed frequency band. The latter estimates the complete differentiated trace series for arbitrary directions in $X$. Actual commutation of the high-regularity flows gives
\[
\mathcal G(\Phi_E^tU_0)-\mathcal G(\Psi_k^tU_0)
=\int_0^tD_U \mathcal G(V(a))[\mathcal V_E-Y_k](V(a))\,da,
\quad V(a)=\Phi_E^a\Psi_k^{t-a}U_0.
\]
We first fix a cutoff in the Green coordinates, then a cutoff in the background fields, and finally let $k$ tend to infinity. The common Fourier tails and the local inverse of $\mathcal G$ give strong convergence in $X$. Comparing two smooth solutions through the same fixed-parameter auxiliary flow shows that their trajectories form a Cauchy sequence in $C([-T,T];X)$. This also gives independence of the approximation.

On the circle, we compute the quadratic trace using the discrete Fourier sums, including the constant spectral mode. The two-branch average again gives a positive endpoint weight, and the full higher-order remainder is integrable in the spectral parameter. Since the fields have zero mean, the inverse from $w$ to $(u,v)$ is bounded. We can therefore combine the spectral bound directly with Kishimoto's local stability estimate and continue on intervals whose length depends only on the endpoint norm. This gives the periodic global flow and its Lipschitz dependence. On both domains, convergence in the local resolution space identifies the differentiated nonlinear source and gives uniqueness in the prescribed integral class.

\subsection{Notation}\label{subsec:notation}

Unless a subscript or Section~\ref{sec:periodic} indicates otherwise, spatial spaces are on $\mathbb R$. We use the unitary Fourier transform on the line and normalized Fourier coefficients on the physical circle. The Schwartz space is denoted by $\mathcal S$ and the space of distributions by $\mathcal D'$. Write
\[
\langle\xi\rangle=(1+\xi^2)^{1/2},\quad
X^s=H^s\times H^{s-1},\quad X=X^{-1/2},\quad
J=\begin{pmatrix}0&\partial_x \\\partial_x &0\end{pmatrix}.
\]
We use $C_tY=C(I;Y)$ on the time interval under consideration and $\tau_yf(x)=f(x+y)$. Multiplication by $f$ is written $M_f$, so that $M_f g=fg$. Spatial derivatives are written $\partial_x$ and $\partial_y$. The notation $D_UF(U)[h]$ means the Fr\'echet derivative of the map $F$ with respect to its input $U$, applied to the direction $h$. For scalar functionals we also use $dA_z(U)[h]$ and the variational derivatives $\delta A_z/\delta u$ and $\delta A_z/\delta v$.
For a complex-valued function $f$, write
\[
\Re f=\tfrac12(f+\bar f),\qquad
\Im f=\tfrac1{2i}(f-\bar f).
\]
The complex variable used in the local theory is
\[
w=u+i(1-\partial_x^2)^{-1}\partial_x v,\qquad
\Re w=u,\quad \Im w=(1-\partial_x^2)^{-1}\partial_x v.
\]
The strong and weak Green spaces are
\[
Z_G=H^{3/2}\times H^{1/2},\qquad W_G=H^{-1/2}\times H^{-3/2}.
\]
We use Kishimoto's endpoint spacetime space $\mathcal W_p$. The subscript $p$ records one fixed choice of the overlapping Fourier projections in its definition.
\begin{definition}\label{def:Wp}
Let $\widetilde b(\tau,\xi)$ be the spacetime Fourier transform, and put
\[
\mu=\langle\tau+\xi^2\rangle,\quad\nu=\langle\xi\rangle,\qquad
\|b\|_{X_{\rm Sch}^{s,c}}=\|\nu^s\mu^c\widetilde b\|_{L^2_{\tau,\xi}},\quad
\|b\|_{Y_{\rm Sch}^s}=\|\nu^s\widetilde b\|_{L^2_\xi L^1_\tau}.
\]
With the fixed overlapping regions
\[
\mathcal R_1=\{\mu\lesssim_p\nu\},\quad
\mathcal R_2=\{\nu\lesssim_p\mu\lesssim_p\nu^2\},\quad
\mathcal R_3=\{\mu\gtrsim_p\nu^2\},\quad
\mathcal R_4=\{\mu\gg_p\nu^2\},
\]
define
\begin{equation}\label{eq:Wp}
\begin{split}
\|b\|_{\mathcal W_p}={}&\|P_{\mathcal R_1}b\|_{X_{\rm Sch}^{-1/2,1}}
+\|P_{\mathcal R_2}b\|_{X_{\rm Sch}^{1/2,0}}\\
&+\sum_{M\ {\rm dyadic}}
\|P_{\mathcal R_3\cap\{\mu\sim_p M\}}b\|_{X_{\rm Sch}^{1/2,0}}
+\|P_{\mathcal R_4}b\|_{Y_{\rm Sch}^{-1/2}}.
\end{split}
\end{equation}
Here $P_{\mathcal R}$ is the fixed spacetime Fourier projection for the indicated region, and the dyadic sum is over $M=1,2,4,\ldots$. For an interval $I$,
\[
\|b\|_{\mathcal W_p(I)}
=\inf\{\|b_{\rm ext}\|_{\mathcal W_p}:b_{\rm ext}|_I=b\}.
\]
On $\mathbb T_\lambda$, use $\xi\in\lambda^{-1}\mathbb Z$, the frequency measure $\lambda^{-1}\sum_\xi$, and
\[
\widehat f(\xi)=(2\pi)^{-1/2}\int_{\mathbb T_\lambda}f(y)e^{-iy\xi}\,dy,
\]
with the unitary time transform. We denote this space by $\mathcal W_{p,\lambda}$. This local Fourier convention differs from the normalized coefficients on the physical circle.
\end{definition}
Restriction is defined by extension, not by applying an ordinary time cutoff. The original vector field and its smooth Hamiltonian are
\[
\mathcal V_E=(\partial_x v,\partial_x (u-\partial_x^2u-u^2)),\qquad
E=\frac12\int(v^2+u^2+u_x^2)\,dx-\frac13\int u^3\,dx.
\]
The linearization at $U=0$ is $\partial_tU=\mathscr L_EU$, with generator and solution group
\begin{equation}\label{eq:linear-flow-notation}
\mathscr L_E=\begin{pmatrix}0&\partial_x \\\partial_x (1-\partial_x^2)&0\end{pmatrix},
\qquad \mathcal U_E(t)=e^{t\mathscr L_E}.
\end{equation}
Thus $\mathcal V_E(U)=\mathscr L_EU+N(U)$, where $N(U)=(0,-\partial_x (u^2))$. The Fourier representation of $\mathcal U_E$ is given in \eqref{eq:legacy-38}.

For the Lax operator, set
\[
\mathfrak a=\tfrac12u-\tfrac14,\qquad \mathfrak b=\tfrac1{4\sqrt3}v,
\qquad L=\partial_x^3+\mathfrak a \partial_x +\tfrac12(\partial_x\mathfrak a)+\mathfrak b=L_0+V_U.
\]
Compositions such as $\partial_x\circ M_f$ act as $g\mapsto\partial_x(fg)$. The functions $\mathfrak a$ and $\mathfrak b$ are the Lax coefficients. We use
\[
d_k=3k^2-\tfrac14,\qquad C_k=8d_k^3/k,\qquad\eta=z_1=\tfrac34.
\]
A spectral point is called admissible if $z=\pm z_k$ for some $k\ge1$. For real $z\ne0$, put $R_{0,z}=(L_0-z)^{-1}$. Its Fourier symbol is
\begin{equation}\label{eq:free-resolvent-symbol}
m_z(p)=\frac{1}{-i(p^3+p/4)-z},\qquad
\widehat{R_{0,z}f}(p)=m_z(p)\widehat f(p),\quad p\in\mathbb R.
\end{equation}
Its convolution kernel is denoted by $g_z$:
\[
g_z(s)=\frac1{2\pi}\int_{\mathbb R}e^{ips}m_z(p)\,dp,
\qquad R_{0,z}f(x)=\int_{\mathbb R}g_z(x-y)f(y)\,dy.
\]
The integral defining $g_z$ is absolutely convergent. The notation $g_z$ is distinct from the data gradients $g_u,g_v$ below.
Choose the continuous square root \(t_z^2=m_z\) with \(\Im t_z>0\) for \(z>0\), and \(\Re t_z>0\) for \(z<0\). Then
\[
\begin{gathered}
\widehat{T_zf}(p)=t_z(p)\widehat f(p),\qquad B_z=T_zV_UT_z,\\
S_z=(I+B_z)^{-1},\qquad R_z=T_zS_zT_z,
\end{gathered}
\]
and, on the common Neumann ball,
\[
A_z=\sum_{n\ge2}\frac{(-1)^{n+1}}{n}\operatorname{Tr}(B_z^n).
\]
Gradients use the real distributional pairing,
\[
dA_z[h]=\langle g_u,h_u\rangle+\langle g_v,h_v\rangle,
\qquad X^*=H^{1/2}\times H^{3/2}.
\]
At the fixed reference point $\eta=z_1=3/4$, define the paired Green map by
\begin{equation}\label{eq:reference-green-map}
\mathcal G(U)=\bigl(4\sqrt3\,g_v(\eta;U),\,4g_u(\eta;U)\bigr)
=\left(4\sqrt3\,\frac{\delta A_\eta}{\delta v}(U),\,
4\frac{\delta A_\eta}{\delta u}(U)\right)\in Z_G.
\end{equation}
Lemma~\ref{lem:determinant} shows that
$\mathcal G:X\to Z_G$ is real analytic on the spectral ball.
Its representation in terms of the free-subtracted Green kernel
is established in
Section~\ref{subsec:green-kernels-and-two-parameter-identities}.
The notation $P_{\le N}$ denotes the orthogonal Fourier projection onto $|\xi|\le N$, and $P_{>N}=I-P_{\le N}$. A nonhomogeneous Littlewood--Paley family is denoted by $P_j$, $j\ge0$, with $P_0$ containing the low frequencies. We write $A\lesssim B$ for $A\le CB$ and $A\sim B$ for two such bounds. Subscripts list additional dependence of the implicit constant; uniformity in $k$ is stated explicitly. Constants may also depend on the fixed local decomposition and reference point.

Section~\ref{sec:preliminaries} gives the local and operator estimates. Section~\ref{sec:spectral} proves the spectral bounds and Green identities, and Section~\ref{sec:flows} constructs the auxiliary flows and establishes commutation. The real-line limit and its nonlinear source are obtained in Section~\ref{sec:green}. Section~\ref{sec:periodic} treats the periodic problem.

\section{Preliminaries}\label{sec:preliminaries}
We recall the linear and bilinear estimates used in Kishimoto's endpoint local well-posedness theorem and prove persistence of higher regularity on a time interval determined by the endpoint norm. We also record the operator-ideal bounds and Sobolev product estimates required for the perturbation determinant.

\subsection{Operator ideals and products}\label{subsec:operator-ideals-and-products}

We record the operator-ideal bounds and Sobolev products needed to define the determinant and its weak variations. The same bounds control all higher trace terms on the small spectral ball.

Let $\mathfrak S_2$ and $\mathfrak S_1$ denote the Hilbert--Schmidt and trace classes; their norms are also written $\|\cdot\|_{\rm HS}$ and $\|\cdot\|_{\mathfrak S_1}$. Bounded operators have norm $\|\cdot\|_{\rm op}$. For $B\in\mathfrak S_2$,
\[
\det{}_2(I+B)=\det\bigl((I+B)e^{-B}\bigr).
\]
We use\[
\|AB\|_{\mathfrak S_1}\le\|A\|_{\mathfrak S_2}\|B\|_{\mathfrak S_2},\qquad
|\operatorname{Tr}(AB)|\le\|A\|_{\mathfrak S_2}\|B\|_{\mathfrak S_2},
\]
and trace cyclicity for trace-class products \cite{r38}.

The following products control the Green formulas, including their negative-order terms. We use a nonhomogeneous paraproduct decomposition as in \cite{r39}.

\begin{lemma}\label{lem:products}
The products initially defined on Schwartz functions extend continuously as follows:
\[
\begin{gathered}
\|fg\|_{H^s}\lesssim\|f\|_{H^{3/2}}\|g\|_{H^s},
\qquad s\in\{-3/2,-1/2,1/2,3/2\},\\
\|fg\|_{H^{-3/2}}\lesssim\|f\|_{H^{1/2}}\|g\|_{H^{-1/2}},\qquad
\|fg\|_{L^2}\lesssim\|f\|_{H^{1/2}}\|g\|_{H^{1/2}}.
\end{gathered}
\]
\end{lemma}
\begin{proof}
Fix a nonhomogeneous Littlewood--Paley decomposition $P_j$. Write $f_j=P_jf$, $f_{<j}=\sum_{a<j}P_af$, and similarly for $g$. Put
\[
F_j=2^{3j/2}\|P_jf\|_2,\qquad G_j=2^{j/2}\|P_jg\|_2.
\]
Low--high and high--low interactions satisfy
\[
\begin{aligned}
\|f_{<j}\|_\infty&\lesssim\|f\|_{H^{3/2}},\\
\|g_{<j}\|_\infty&\lesssim\sqrt{j+1}\|g\|_{H^{1/2}},\\
2^{j/2}\|f_jg_{<j}\|_2&\lesssim
2^{-j}\sqrt{j+1}F_j\|g\|_{H^{1/2}}.
\end{aligned}
\]
For \(|j-j'|\lesssim1\), \(\ell\le j+O(1)\), Bernstein gives
\[
2^{\ell/2}\|P_\ell(f_jg_{j'})\|_2
\lesssim2^{\ell-2j}F_jG_{j'}.
\]
The kernel \(2^{-(j-\ell)}2^{-j}\) is summable. Discrete Young and Cauchy--Schwarz imply
\[
\|fg\|_{H^{1/2}}\lesssim\|f\|_{H^{3/2}}\|g\|_{H^{1/2}}.
\]
The same decomposition gives the \(H^{3/2}\) algebra estimate, with high--high weight \(2^{2\ell-3j}\). By duality and \(H^{1/2}\hookrightarrow L^4\),
\[
\begin{gathered}
H^{3/2}\cdot H^{-1/2}\subset H^{-1/2},\qquad
H^{3/2}\cdot H^{-3/2}\subset H^{-3/2},\\
H^{1/2}\cdot H^{-1/2}\subset H^{-3/2},\qquad
H^{1/2}\cdot H^{1/2}\subset L^2.
\end{gathered}
\]
Smooth density gives the stated continuous extensions. Adding a constant to the $H^{3/2}$ multiplier is also allowed. This completes the proof.
\end{proof}

\subsection{Endpoint local solutions and persistence}\label{subsec:local}
We use Kishimoto's linear and bilinear estimates for the real-line problem \cite[Section 2]{r1}. Throughout this subsection, $W=\mathcal W_p$ is the resolution space of Definition~\ref{def:Wp}. In the rescaled equation for $b$, the spatial variable is $y$ and $\xi$ is its Fourier dual variable. We recall the Duhamel extension and the estimates needed to identify the limit of smooth solutions and to prove persistence of regularity. Set
\[
\alpha=\tau+\xi^2,\qquad
\Omega_\lambda(\xi)=\frac{\lambda^2\xi^2}{1+\lambda^2\xi^2}.
\]
Fix \(\psi\in C_c^\infty(\mathbb R)\), equal to one on \([-1,1]\) and supported in \([-2,2]\). The same function is used for the time and modulation cutoffs. Define the Duhamel extension by
\[
\begin{aligned}
\mathcal IF={}&-\frac{i}{\sqrt{2\pi}}\psi(t)e^{it\partial_y^2}\mathcal F_\xi^{-1}
\int\widetilde F\psi(\alpha)\sum_{n\ge1}\frac{t^n(i\alpha)^{n-1}}{n!}\,d\tau\\
&+\frac{i}{\sqrt{2\pi}}\psi(t)e^{it\partial_y^2}\mathcal F_\xi^{-1}
\int\widetilde F\frac{1-\psi(\alpha)}{i\alpha}\,d\tau
-i\mathcal F_{\tau,\xi}^{-1}\!\left[\widetilde F\frac{1-\psi(\alpha)}{i\alpha}\right].
\end{aligned}
\]
For \(|t|<1\), it equals \(-i\int_0^te^{i(t-a)\partial_y^2}F(a)\,da\). We use the linear estimates, the bound for this Duhamel extension, the short-interval estimate for solutions with zero initial data, and the following bilinear estimates:
\[
\begin{aligned}
\|\Lambda^{-1}\Omega_\lambda(b\bar c)\|_W
&\le C_p\lambda^{1/2}\|b\|_W\|c\|_W,\\
\|\Lambda^{-1}\Omega_\lambda(bc)\|_W+
\|\Lambda^{-1}\Omega_\lambda(\bar b\bar c)\|_W
&\le C_p\|b\|_W\|c\|_W,
\qquad\Lambda=\langle\tau+\xi^2\rangle.
\end{aligned}
\]
These estimates are proved in \cite[Section 2]{r1}.

For the real-line extension equation at scale $\lambda$, we likewise write
\[
\begin{aligned}
A_{\rm loc}f&=\psi(t)e^{it\partial_y^2}f,\\
L_{\rm loc}b&=\tfrac12\lambda^{-2}\mathcal I(b-\bar b),\\
B_{\rm loc}(b,c)&=-\mathcal I\bigl[\Omega_\lambda(\Re b\,\Re c)\bigr].
\end{aligned}
\]
The bilinear map is first defined on smooth functions and then extended by the preceding estimates.

For an interval $I$, the class $\mathcal K_{\lambda_*}(I)$ consists of real trajectories $U\in C(I;X)$ for which, near every $t_0\in I$,
\[
b(s,y)=\lambda_*^{-2}w(t_0+\lambda_*^{-2}s,y/\lambda_*)
\]
belongs to the restriction space $W$ and satisfies the Duhamel equation associated with $A_{\rm loc}$, $L_{\rm loc}$, and $B_{\rm loc}$ on that window. The scale $\lambda_*$ is fixed for the whole trajectory. We do not assume boundedness of dilation on the full-time space $W$.

The contraction must preserve high regularity without a lifespan depending on high norms. The next lemma obtains this from the analytic fixed-point map and then recovers the second field.

\begin{lemma}\label{lem:persistence}
Fix $\lambda_*$ sufficiently large and then $\varepsilon$ sufficiently small. There is a fixed extension map for the complex local equation on $[-\tau_0,\tau_0]$, $\tau_0=\lambda_*^{-2}$, which is real analytic into $\mathcal W_p$. For integers $n\ge0$,
\[
w_0\in H^{n-1/2},\quad\|U_0\|_X\le2\varepsilon
\quad\Longrightarrow\quad
w\in C([-\tau_0,\tau_0];H^{n-1/2}).
\]
The same $\tau_0$ is used at every higher order. For $m=7/2,11/2$, the recovered real two-field solution satisfies
\[
U_0\in X^m\quad\Longrightarrow\quad
U\in C_tX^m\cap C_t^1X^{m-2}
\]
on this interval, and agrees with the high-regularity mild solution.
\end{lemma}

\begin{proof}
We first construct the analytic extension map on a fixed endpoint lifespan. We then prove persistence of the complex solution and recover the two-field regularity and mild uniqueness.

For \(w=u+i(1-\partial_x^2)^{-1}\partial_x v\) and \(b_0(y)=\lambda^{-2}w_0(y/\lambda)\),
\[
\|w_0\|_{H^{-1/2}}\le\|U_0\|_X,\qquad
\|b_0\|_{H^{-1/2}}\le\lambda^{-1}\|U_0\|_X.
\]

The integral extension gives the equation on the real Banach space \(\mathcal W_p\):\begin{equation}\label{eq:legacy-2}
\begin{gathered}
b=A_{\rm loc}b_0+L_{\rm loc}b+B_{\rm loc}(b,b),\\
\|A_{\rm loc}\|\le C_p,\ \|L_{\rm loc}\|\le C_p\lambda^{-2},\
\|B_{\rm loc}\|\le C_p\lambda^{1/2}.
\end{gathered}
\end{equation}
The quadratic map is real bilinear. Choose
\[
C_p\lambda^{-2}\le\tfrac18,\qquad
8C_p^2\varepsilon\lambda^{-1/2}\le\tfrac18,
\qquad R_W=4C_p\varepsilon/\lambda.
\]
For \(\|U_0\|_X\le2\varepsilon\), the fixed-point map preserves \(B_W(0,R_W)\) and
\[
\|L_{\mathrm{loc}}+2B_{\mathrm{loc}}(b,\cdot)\|\le q\le\tfrac14,
\qquad \|(I-L_{\mathrm{loc}}-2B_{\mathrm{loc}}(b,\cdot))^{-1}\|\le(1-q)^{-1}.
\]
The implicit function theorem gives a real-analytic extension map on the physical interval \(\tau_0=\lambda^{-2}\).

Translations are strongly continuous on each component of \(W\), including the dyadic \(\ell^1\) sum. The extension map \(\mathcal M\) is equivariant, and
\[
W\hookrightarrow Y_{\mathrm{Sch}}^{-1/2}\hookrightarrow C_tH^{-1/2},
\qquad \mathcal M(\tau_yb_0)=\tau_y\mathcal M(b_0).
\]
For \(b_0\in H^{n-1/2}\), differentiate the translation orbit in \(H^{-1/2}\). Distributional difference quotients identify
\[
\partial_y^n b=\left.\partial_h^n\mathcal M(\tau_hb_0)\right|_{h=0}\in W.
\]

At the highest order, \(L_{\rm loc}\partial_y^n b+2B_{\rm loc}(b,\partial_y^n b)\) is absorbed using the same contraction constant $q$:\begin{equation}\label{eq:legacy-3}
\|\partial_y^nb\|_W\le(1-q)^{-1}
[C_p\|\partial_y^nb_0\|_{H^{-1/2}}+
C_p\lambda^{1/2}\sum_{j=1}^{n-1}\binom nj
\|\partial_y^jb\|_W\|\partial_y^{n-j}b\|_W].
\end{equation}
Induction in \eqref{eq:legacy-3} gives
\[
w\in C([-\tau_0,\tau_0];H^{n-1/2}),
\]
and the same argument for the joint directional maps gives \(C^2\) dependence in higher Sobolev spaces. The interval \(\tau_0\) is independent of their norms. For \(m=7/2,11/2\), approximation in \(H^m\) and the algebra estimate imply
\[
u_j\to u\text{ in }C_tH^m
\quad\Longrightarrow\quad
u_j^2\to u^2\text{ in }C_tH^m,
\]
identifying the high-regularity local source with the ordinary product.

Define\begin{equation}\label{eq:legacy-4}
v(t)=v_0+\int_0^t\partial_x (u-\partial_x^2u-u^2)\,da
\quad\text{initially in }C_tH^{m-3}.
\end{equation}
The complex equation and the initial compatibility give
\[
u_t=(1-\partial_x^2)\Im w,
\qquad u_{tt}=\partial_x^2(u-\partial_x^2u-u^2),\qquad \partial_x v=u_t.
\]

At high frequencies,\[
P_{>1}v=P_{>1}\partial_x^{-1}(1-\partial_x^2)\Im w\in C_tH^{m-1}.
\]
Here \(P_{>1}\partial_x^{-1}(1-\partial_x^2)\) has order one.

At low frequencies, estimate \eqref{eq:legacy-4} directly:\begin{equation}\label{eq:legacy-5}
\|P_{\le1}v(t)\|_{H^{m-1}}\le\|P_{\le1}v_0\|_{H^{m-1}}+
C_m|t|\sup_a(\|u(a)\|_{H^m}+\|u(a)\|_{H^m}^2).
\end{equation}
Consequently,
\[
U\in C_tX^m\cap C_t^1X^{m-2},\qquad
\operatorname{Dom}(\mathscr L_E:X^{m-2}\to X^{m-2})=X^m.
\]
Differentiating \(\mathcal U_E(-t)U(t)\) gives
\[
U(t)=\mathcal U_E(t)U_0+\int_0^t\mathcal U_E(t-a)(0,-\partial_x (u(a)^2))\,da.
\]
For \(m>1/2\), the nonlinearity is locally Lipschitz \(X^m\to X^m\); Gronwall identifies the high mild solution on the entire low-norm interval.
The fixed lifespan, persistence, and high-regularity identification have now been proved. This completes the proof.
\end{proof}

\begin{remark}\label{rem:persistence-use}
The cases $m=7/2$ and $m=11/2$ in Lemma~\ref{lem:persistence} provide the regularity used for spectral conservation and the common-trajectory chain rule, respectively. Section~\ref{sec:flows} develops the corresponding high-space variations. For the general smooth-approximation method, see \cite{r40}.
\end{remark}

\section{Spectral conservation and two-parameter identities}\label{sec:spectral}
We define the balanced determinant at the endpoint, prove coercivity of the two-branch average, and derive conserved tails. The matrix-kernel construction then justifies the two-parameter identities used in Sections~\ref{sec:flows} and~\ref{sec:green}.

\subsection{Balanced resolvents and determinant variations}\label{subsec:balanced-resolvents-and-determinant-variations}

We first balance the third-order free inverse on both sides of the perturbation. The resulting Hilbert--Schmidt estimate is valid for endpoint coefficients and supplies the small parameter for the full trace series.

\begin{lemma}\label{lem:free-resolvent}
Let $k\ge1$ and $z=\pm z_k$. The free resolvent $R_{0,z}$ has Fourier symbol
\[
m_z(p)=\frac{1}{-i(p^3+p/4)-z},\qquad p\in\mathbb R,
\]
as in \eqref{eq:free-resolvent-symbol}. Its convolution kernel exists, $|m_z(p)|\sim(k^3+|p|^3)^{-1}$, and\begin{equation}\label{eq:legacy-6}
\|B_z(U)\|_{\rm HS}^2\le K^2\int k^{-2}[(k+|\xi|)^{-1}|\widehat u|^2+(k+|\xi|)^{-3}|\widehat v|^2]d\xi
\le K^2 k^{-2}\|U\|_X^2.
\end{equation}
The constant $K$ is absolute and independent of the data, scale, and spectral sign.
\end{lemma}

\begin{proof}
We construct the free convolution kernel and establish the multiplier comparison, then prove the Hilbert--Schmidt estimate \eqref{eq:legacy-6}.

For \(f(r)=r^3-r/4\), the roots of \(f(r)=z_k\) are
\[
r_1=k,\qquad r_{2,3}=-k/2\pm\tfrac i2\sqrt{3k^2-1}.
\]
Partial fractions give
\[
g_z(x)=\begin{cases}
\displaystyle\sum_{\Re r_j<0}\frac{e^{r_jx}}{f'(r_j)},&x>0,\\
\displaystyle-\sum_{\Re r_j>0}\frac{e^{r_jx}}{f'(r_j)},&x<0,
\end{cases}
\qquad g_{-z}(x)=-g_z(-x).
\]
For a function $q$ with one-sided limits at zero, $[q]=q(0+)-q(0-)$. Since
\[
\sum_j\frac1{f'(r_j)}=\sum_j\frac{r_j}{f'(r_j)}=0,
\qquad\sum_j\frac{r_j^2}{f'(r_j)}=1,
\]
we have
\[
[g]=[g']=0,\qquad[g'']=1,\qquad(L_0-z)g_z=\delta,
\qquad|g_z^{(j)}(x)|\lesssim k^{j-2}e^{-k|x|/2}\quad(j\le2).
\]
Also \(3k^3/4\le z_k\le k^3\), giving the multiplier comparison.

The kernel is\[
B_z(p,q)=(2\pi)^{-1/2}t_z(p)t_z(q)[i(p+q)\widehat u(p-q)/4+\widehat v(p-q)/(4\sqrt3)].
\]
Put \(R=|\xi|\) and \(\omega_k(p)=(k^3+|p|^3)^{-1}\). Scaling for \(R\le k\), and the decomposition
\[
\{|p|\le R/2\}\cup\{|p-\xi|\le R/2\}\cup
\{|p|,|p-\xi|>R/2\}
\]
for \(R>k\), yield
\[
\begin{aligned}
\int\omega_k(p)\omega_k(p-\xi)\,dp
&\lesssim k^{-2}(k+R)^{-3},\\
\int(2p-\xi)^2\omega_k(p)\omega_k(p-\xi)\,dp
&\lesssim k^{-2}(k+R)^{-1}.
\end{aligned}
\]
Indeed, the low regions use \(\int\omega_k\lesssim k^{-2}\); the remaining region scales as \(R^{-5}\) and \(R^{-3}\). The vertex bound \(|a+b|^2\le2|a|^2+2|b|^2\) proves \eqref{eq:legacy-6}. Approximation of the Fourier kernel in \(L^2(\mathbb R^2)\) defines \(B_z\) for every \(U\in X\).
The kernel construction, multiplier comparison, and \eqref{eq:legacy-6} prove the lemma.
\end{proof}

The determinant must be differentiated before any classical rough kernel diagonal is invoked. The next statement gives its derivatives directly in the operator ideals and identifies their distributional dual spaces.

\begin{lemma}\label{lem:determinant}
Fix \(\beta _*<1\) and take \(2K\varepsilon \le \beta _*\). On the ball of radius \(2\varepsilon \), the functional \(A_z\) is real analytic and $\det_2$ is nonzero. For every \(h,l\in X\),\begin{equation}\label{eq:legacy-7}
dA_z[h]=-\operatorname{Tr}(S_zB_zB_z(h)),\quad
d^2A_z[h,l]=-\operatorname{Tr}(S_zB_z(l)S_zB_z(h)),
\end{equation}\begin{equation}\label{eq:legacy-8}
|dA_z[h]|\le C_{\beta_*} k^{-2}\|U\|\|h\|,
\quad |d^2A_z[h,l]|\le C_{\beta_*} k^{-2}\|h\|\|l\|.
\end{equation}
In particular, \(g_u\in H^{1/2}\) and \(g_v\in H^{3/2}\), with a Lipschitz bound uniform in \(k^{-2}\) on the ball.
\end{lemma}

\begin{proof}
We prove analyticity and nonvanishing of the determinant, derive \eqref{eq:legacy-7}--\eqref{eq:legacy-8}, and then verify reality and the asserted gradient regularity.

Suppress the spectral subscript on $B_z$ and $S_z$. The operator $ (I+B)e^{-B}-I$ is trace class. Finite-rank approximation and the Neumann inverse identify the determinant branch. Moreover,
\[
|\operatorname{Tr}B^n|\le\|B\|_{\mathfrak S_2}^2\beta_*^{n-2},
\qquad\|S\|\le(1-\beta_*)^{-1},\qquad dS[h]=-SB(h)S.
\]
Every differentiated trace contains two Hilbert--Schmidt factors. Geometric summation in \(\mathfrak S_1\) therefore permits differentiation and gives \eqref{eq:legacy-7}--\eqref{eq:legacy-8}, together with all higher variations. The gradients belong to the distributional dual \(X^*\).

For real data,
\[
t_z(-p)=\begin{cases}-\overline{t_z(p)},&z>0,\\\overline{t_z(p)},&z<0.
\end{cases}
\]
Thus \(B_z\) commutes with \(f(p)\mapsto\overline{f(-p)}\), and
\[
\operatorname{Tr}B_z^n\in\mathbb R,\qquad A_z\in\mathbb R.
\]
Thus \eqref{eq:legacy-7}--\eqref{eq:legacy-8} and all the asserted determinant and gradient properties hold. This completes the proof.
\end{proof}

\begin{remark}\label{rem:weak-realization}
The balanced factorization also gives
\[
L-z=T_z^{-1}(I+B_z)T_z^{-1}:H^{3/2}\to H^{-3/2},
\qquad (L-z)^{-1}=R_z.
\]
This weak realization agrees with the differential operator for smooth coefficients.
\end{remark}

\subsection{Coercivity, tails, and conservation}\label{subsec:coercivity-tails-and-conservation}

A single spectral branch is not coercive. We therefore compute the quadratic term before estimating the higher-order series; the branch average supplies the positive weight needed for the endpoint barrier.

\begin{lemma}\label{lem:coercivity}
Set \(\chi_k(\xi)=3k^2+1/4+\xi^2\) and \(\mathfrak q_k(\xi)=\xi^4+(3k^2+1/2)\xi^2+d_k^2\). The quadratic term is exactly\begin{equation}\label{eq:legacy-9}
A_{z_k}^{[2]}=-\int\left[\frac{k\langle\xi\rangle^2|\widehat u|^2+k|\widehat v|^2}{16d_k\mathfrak q_k}
+\frac{\chi_k\Re(\widehat u\overline{\widehat v})}{8\sqrt3d_k\mathfrak q_k}\right]d\xi.
\end{equation}
Changing the spectral sign reverses only the mixed term. Thus \(\mathscr S^{[2]}=\int h(\xi)(|\widehat u|^2+|\widehat v|^2/\langle\xi\rangle^2)d\xi\), with \(h\sim\langle\xi\rangle^{-1}\). For small \(\varepsilon \), there are \(\gamma ,\Gamma >0\) such that\begin{equation}\label{eq:legacy-10}
\gamma\|U\|_X^2\le\mathscr S(U)\le\Gamma\|U\|_X^2,
\quad \|P_{>R}U\|_X^2\le C(|\mathscr S_R(U)|+\varepsilon^3/R).
\end{equation}
The constants are absolute or depend on \(\beta _*\) and the chosen \(\varepsilon \). Here \(P_{\le N}\) is the orthogonal Fourier projection onto \(|\xi |\le N\).
\end{lemma}

\begin{proof}
We compute \eqref{eq:legacy-9} and the averaged quadratic weight, then bound the complete trace remainder to obtain both inequalities in \eqref{eq:legacy-10}.

In this proof, write $A_\pm=A_{\pm z_k}$ and $m_\pm=m_{\pm z_k}$.
Denote the balanced Hilbert--Schmidt constant by \(K\ge 1\); then \(\|B_{\pm z_k}(U)\|_{HS}\le Kk^{-1}\|U\|_X\).
Fix \(\beta_*<1\) with \(\|B_z(U)\|_{\mathrm{op}}\le\beta_*\) on a slightly larger spectral ball.

We compute the quadratic trace using the three integrals \(I_j=(2\pi)^{-1}\int(2p+\xi)^j m_+(p)m_+(p+\xi)dp\).
The integrand \(j\le 2\) is absolutely integrable at infinity. The two lower half-plane poles give\[
I_0=\frac{6k}{d_k\mathfrak q_k},\quad I_1=-\frac{2i(3k^2+1/4+\xi^2)}{d_k\mathfrak q_k},
\quad I_2=-\frac{2k(1+\xi^2)}{d_k\mathfrak q_k}.
\]
Multiply by the vertex \(i(p+q)\widehat u/4+\widehat v/(4\sqrt3)\) and take \(-\text{Tr}(B^2)/2\). The positive branch gives\begin{equation}\label{eq:legacy-11}
A_+^{[2]}=-\int\frac{k(\langle\xi\rangle^2|\widehat u|^2+|\widehat v|^2)}
{16d_k\mathfrak q_k}d\xi
-\int\frac{(3k^2+1/4+\xi^2)\Re(\widehat u\overline{\widehat v})}
{8\sqrt3d_k\mathfrak q_k}d\xi .
\end{equation}
The negative branch reverses only the mixed term. Hilbert--Schmidt continuity gives reality and extension of the integral to rough data. The individual branches need not be positive.
Set \(s_k=-(A_++A_-)/2\). Then\begin{equation}\label{eq:legacy-12}
\begin{gathered}
s_k^{[2]}=\int c_k(\xi)(|\widehat u|^2+|\widehat v|^2/\langle\xi\rangle^2)d\xi,
\\ c_k=\frac{k\langle\xi\rangle^2}{16d_k\mathfrak q_k}.
\end{gathered}
\end{equation}
For \(k\ge 1\),\[
\tfrac{11}4k^2\le d_k\le3k^2,\quad
(k^2+\xi^2)^2\le \mathfrak q_k\le9(k^2+\xi^2)^2.
\]
Consequently,\begin{equation}\label{eq:legacy-13}
\frac{\langle\xi\rangle^2}{432(k^2+\xi^2)^2}
\le k c_k(\xi)\le
\frac{\langle\xi\rangle^2}{44(k^2+\xi^2)^2}.
\end{equation}
The full integral starts at $1$ and the tail starts at \(R\ge 1\).
We have \(h_R(\xi)=\int_R^\infty k c_k(\xi)dk\). Split into \(|\xi|\le 1\) and $>1$. In the latter case, integrate over \([|\xi|,2|\xi|]\) for a lower bound. Alternatively, take \(a=\max(1,|\xi|)\); on $[a,2a]$, \(k^2+\xi^2\le5a^2\). Thus\begin{equation}\label{eq:legacy-14}
c_2\langle\xi\rangle^{-1}\le h_1(\xi)\le C_2\langle\xi\rangle^{-1},
\quad c_2=1/10800,\quad C_2=1.
\end{equation}
For the upper bound, use \(\int_1^\infty k^{-4}\,dk=1/3\) when \(|\xi|\le1\), and scale \(k=|\xi|x\) otherwise. For \(|\xi|>R\), integration over \([|\xi|,2|\xi|]\) gives the same lower bound for \(h_R\).

Moreover, for every \(A\ge 1\),\begin{equation}\label{eq:legacy-15}
|\xi|\le A\ \Longrightarrow\ h_R(\xi)\le\langle\xi\rangle^2/(132R^3).
\end{equation}

\textit{The full trace remainder.} Bound \(\operatorname{Tr}B^n\) using two Hilbert--Schmidt factors and operator norms for the remaining factors. For \(n\ge 3\),\begin{equation}\label{eq:legacy-16}
|A_z-A_z^{[2]}|\le\frac{\|B\|_{HS}^2\|B\|_{op}}{3(1-\|B\|_{op})}
\le C_3 k^{-3}\|U\|_X^3,\quad C_3=\frac{K^3}{3(1-\beta_*)}.
\end{equation}
The averaged branches satisfy the same bound. In the derivative series, the factor $n$ cancels $1/n$, and geometric summation gives\begin{equation}\label{eq:legacy-17}
|d(s_k-s_k^{[2]})(U)[h]|\le3C_3 k^{-3}\|U\|_X^2\|h\|_X.
\end{equation}
These estimates sum all \(n\ge 3\) in trace norm, including the differentiated series. They dominate difference quotients on the spectral ball and justify differentiation of the full remainder.
The quantity \(\mathscr S_R=\int_R^\infty ks_kdk\) is absolutely defined, since\begin{equation}\label{eq:legacy-18}
\int_R^\infty k|s_k(U)|dk\le
\mathscr S_R^{[2]}(U)+C_3\|U\|_X^3/R<\infty
\end{equation}
The exact averaged quadratic weight is needed here. Integrating the crude single-branch Hilbert--Schmidt bound \(k^{-2}\) against $k$ would not suffice.
Equations \eqref{eq:legacy-14} and \eqref{eq:legacy-16} give\begin{equation}\label{eq:legacy-19}
|\mathscr S-\mathscr S^{[2]}|\le C_3\|U\|_X^3,\quad
|\mathscr S_R-\mathscr S_R^{[2]}|\le C_3\|U\|_X^3/R .
\end{equation}
Choose \(\varepsilon \) to satisfy \(2K\varepsilon\le\beta_*,\ C_3\varepsilon\le c_2/2\), within the radii allowed by the reference inverse and the local contraction.
The constants can then be chosen as\begin{equation}\label{eq:legacy-20}
\gamma=c_2/2,\quad\Gamma=C_2+C_3\varepsilon,\quad
\gamma\|U\|_X^2\le\mathscr S(U)\le\Gamma\|U\|_X^2,
\end{equation}\begin{equation}\label{eq:legacy-21}
\|P_{>R}U\|_X^2\le c_2^{-1}(|\mathscr S_R(U)|+C_3\varepsilon^3/R).
\end{equation}
Equations \eqref{eq:legacy-11}--\eqref{eq:legacy-14} give the stated quadratic formula and weight, while \eqref{eq:legacy-20}--\eqref{eq:legacy-21} give \eqref{eq:legacy-10}. This completes the proof.
\end{proof}

The preceding estimate converts spectral-tail control into Fourier-tail control. For the approximation argument, we also need continuity of the spectral tails and their uniform decay on precompact sets of initial data.

\begin{lemma}\label{lem:spectral-tail-continuity}
Let $\varepsilon$, $C_2$, and $C_3$ be as in Lemma~\ref{lem:coercivity}. For $R\ge1$ and $U,V\in X$ with $\|U\|_X,\|V\|_X\le\varepsilon$,
\[
|\mathscr S_R(U)-\mathscr S_R(V)|
\le \left(2C_2\varepsilon+\frac{3C_3\varepsilon^2}{R}\right)\|U-V\|_X.
\]
If $Q\subset\{U\in X:\|U\|_X\le\varepsilon\}$ is precompact and
\[
\rho_Q(A)=\sup_{U\in Q}\|P_{>A}U\|_X,
\]
then, for $A,R\ge1$,
\begin{equation}\label{eq:legacy-22}
\sup_{U\in Q}|\mathscr S_R(U)|\le C_2\rho_Q(A)^2+
\frac{\langle A\rangle^3\varepsilon^2}{132R^3}
+\frac{C_3\varepsilon^3}{R}.
\end{equation}
In particular,
\[
\lim_{R\to\infty}\sup_{U\in Q}|\mathscr S_R(U)|=0.
\]
\end{lemma}

\begin{proof}
We prove the Lipschitz estimate first, then derive \eqref{eq:legacy-22} and uniform decay on precompact sets.

Since $0\le h_R\le h_1\le C_2\langle\xi\rangle^{-1}$, the quadratic term satisfies
\[
|\mathscr S_R^{[2]}(U)-\mathscr S_R^{[2]}(V)|
\le C_2(\|U\|_X+\|V\|_X)\|U-V\|_X.
\]
Integrating \eqref{eq:legacy-17} in $k$ and along the segment from $V$ to $U$ gives
\[
\begin{aligned}
&|(\mathscr S_R-\mathscr S_R^{[2]})(U)
  -(\mathscr S_R-\mathscr S_R^{[2]})(V)|\\
&\qquad\le
3C_3\|U-V\|_X
\int_R^\infty k^{-2}\,dk
\int_0^1\|V+\theta(U-V)\|_X^2\,d\theta
\le\frac{3C_3\varepsilon^2}{R}\|U-V\|_X.
\end{aligned}
\]
This proves the Lipschitz estimate.

A finite net of the compact $X$ closure of $Q$, together with the contraction property of $P_{>A}$, gives $\rho_Q(A)\to0$ as $A\to\infty$.
Splitting the quadratic integral at $|\xi|=A$ and using \eqref{eq:legacy-14}--\eqref{eq:legacy-15}, we obtain
\[
\begin{aligned}
\mathscr S_R^{[2]}(U)
&\le C_2\|P_{>A}U\|_X^2
+\frac{1}{132R^3}\int_{|\xi|\le A}
\langle\xi\rangle^2\left(|\widehat u|^2+
\frac{|\widehat v|^2}{\langle\xi\rangle^2}\right)\,d\xi\\
&\le C_2\rho_Q(A)^2+
\frac{\langle A\rangle^3\varepsilon^2}{132R^3},
\qquad U\in Q.
\end{aligned}
\]
Estimate \eqref{eq:legacy-19} yields \eqref{eq:legacy-22}. First choose $A$ large, then let $R\to\infty$.
This proves the stated continuity and uniform decay.
\end{proof}

\begin{remark}\label{rem:spectral-tail-bounded-sets}
The uniform decay in Lemma~\ref{lem:spectral-tail-continuity} need not hold on a bounded small-data set. Choose symmetric real frequency packets with \(\widehat u_j=\delta\langle\xi\rangle^{1/2}
[\psi(\xi-2^j)+\psi(\xi+2^j)]/(\sqrt2\|\psi\|_2)\),v=0,
where \(\psi \) is real and even, supported in \((-1/8,1/8)\). The $X$ norm is exactly \(\delta \). The packets are orthogonal in $X$, so the sequence is not precompact.
Choose \(\delta \) small enough that all data lie in $\mathcal D$. For a fixed large $R$, take \(2^j-1/8>R\).
The lower tail bound in \eqref{eq:legacy-14} and estimate \eqref{eq:legacy-19} give \(\mathscr S_R(U_j)\ge c_2\delta^2-C_3\delta^3/R\).
Hence a bounded small-data set need not have a uniform spectral tail. This does not contradict \eqref{eq:legacy-22}.
\end{remark}

The following conservation statement uses only an existing smooth trajectory. This order of proof will allow continuation without assuming a rough global flow.

\begin{lemma}\label{lem:conservation}
Let $I$ be a compact interval and let $U$ be a smooth solution such that
\[
U\in C(I;H^3\times H^2)\cap C^1(I;H^1\times L^2),\qquad
\sup_{t\in I}\|B_{\pm z_k}(U(t))\|_{\rm op}\le\beta<1.
\]
For every fixed $k\ge1$,
\[
A_{\pm z_k}(U(\cdot))\in C^1(I),\qquad
\frac{d}{dt}A_{\pm z_k}(U(t))=0.
\]
No spatial weights are required.
\end{lemma}

\begin{proof}
For each fixed spectral point, we show that $A_z(U(t))$ is continuously differentiable and that its derivative vanishes. The trace calculation is justified before passing to the infinite series.

Fix $z$ and suppress its subscript on $T_z$, $B_z$, and the resolvents. Set $\mathcal P_E=\sqrt3\partial_x^2+(u-1/2)/\sqrt3$. Direct expansion gives
\[
L_t=[\mathcal P_E,L],\qquad
u_t=v',\qquad v_t=(u-u''-u^2)'.
\]
To justify the trace, write \(a_0(p)=-i(p^3+p/4)\) and
\[
\mathcal C=T[M_{u/\sqrt3},L_0]T,
\qquad \mathcal H_+=TM_{u/\sqrt3}T^{-1},\quad
\mathcal H_-=T^{-1}M_{u/\sqrt3}T.
\]
The last two operators are defined by their Fourier kernels. Schur's test and the three-region decomposition give
\[
\begin{gathered}
|t(p)/t(q)|\lesssim_z\langle p-q\rangle^{3/2},\qquad
\|\mathcal H_\pm\|_{\mathrm{op}}\lesssim_z\|u\|_{H^3},\\
\int|m(p)m(p-\xi)|\,|a_0(p)-a_0(p-\xi)|^2\,dp
\lesssim_z\langle\xi\rangle^3,\\
\|\mathcal C\|_{\mathfrak S_2}+\|[\partial_x^2,B]\|_{\mathfrak S_2}
\lesssim_z\|u\|_{H^2}+\|v\|_{H^1},\qquad
\mathcal C=\mathcal H_+-\mathcal H_-.
\end{gathered}
\]
The corresponding vertex integrals with \(\xi^2(2p-\xi)^4\) and \(\xi^2(2p-\xi)^2\) have weights \(\langle\xi\rangle^3\) and \(\langle\xi\rangle\). Hence \(B_t\in C_t\mathfrak S_2\), and, with \(\mathcal P_0=\sqrt3\partial_x^2\),
\begin{equation}\label{eq:legacy-23}
B_t=[\mathcal P_0,B]+[\mathcal H_+,B]+(I+B)\mathcal C.
\end{equation}
High Sobolev density extends this identity from compact smooth coefficients.

For \(E_\delta=e^{i\delta\mathcal P_0}\), dominated convergence in \(\mathfrak S_2\) gives
\[
\delta^{-1}(E_\delta BE_\delta^{-1}-B)\to i[\mathcal P_0,B],
\qquad\operatorname{Tr}(E_\delta BE_\delta^{-1})^n=\operatorname{Tr}B^n.
\]
Thus \(\operatorname{Tr}(B^{n-1}[\mathcal P_0,B])=0\) for \(n\ge2\). The bounded \(\mathcal H_+\) commutator vanishes by cyclicity. For the partial logarithmic sum,
\[
\partial_t\sum_{n=2}^N\frac{(-1)^{n+1}}n\operatorname{Tr}B^n
=-\operatorname{Tr}(B\mathcal C)+(-1)^{N+1}\operatorname{Tr}(B^N\mathcal C),
\]
\[
|\operatorname{Tr}(B^N\mathcal C)|
\le\|B\|_{\mathfrak S_2}\|\mathcal C\|_{\mathfrak S_2}\beta_*^{N-1}.
\]
The series and its derivative converge uniformly on compact time intervals. Finally,
\[
m(p)m(q)(a_0(p)-a_0(q))=m(q)-m(p),
\]
\[
\begin{aligned}
\int[m(q)-m(q+\xi)]\,dq&=0,\\
\int(2q+\xi)[m(q)-m(q+\xi)]\,dq&=2\xi\int m(q)\,dq.
\end{aligned}
\]
Since \(m,qm\in L^1\), the changes of variables are legitimate. The remaining trace is a constant multiple of
\[
\int i\xi\widehat u(\xi)\widehat u(-\xi)\,d\xi=0.
\]
Therefore $A_z\in C^1(I)$ and $\partial_tA_z=0$ for each fixed admissible spectral point.
This proves the lemma.
\end{proof}

\begin{remark}\label{rem:integrated-spectral-conservation}
If the trajectory remains in the common spectral ball, the pointwise equalities from Lemma~\ref{lem:conservation} and the absolute convergence in \eqref{eq:legacy-18} imply
\[
\mathscr S_R(U(t))=\int_R^\infty ks_k(U(t))\,dk
=\int_R^\infty ks_k(U(0))\,dk=\mathscr S_R(U(0)).
\]
Taking $R=1$ also gives conservation of $\mathscr S$. No differentiation under the spectral integral is needed.
\end{remark}

\subsection{Green kernels and two-parameter identities}\label{subsec:green-kernels-and-two-parameter-identities}

We now connect the abstract trace variations with Green kernels on smooth backgrounds. The diagonal regularity is proved first, so the subsequent two-parameter formulas have well-defined operator domains and diagonal values.

Throughout the kernel calculations, $\mathfrak a$ and $\mathfrak b$ denote the Lax coefficients defined in Section~\ref{subsec:notation}.

For smooth data, let $G_z$ be the kernel of $R_z$ and set
\[
\Delta G_z(x,y)=G_z(x,y)-g_z(x-y),
\]
where $g_z$ is the free resolvent kernel defined in Section~\ref{subsec:notation} and computed in Lemma~\ref{lem:free-resolvent}.
For a kernel $F$ with one-sided limits at the diagonal, we use
\[
F|_{\mathrm{diag},\pm}(y)=\lim_{h\downarrow0}F(y\pm h,y),
\qquad [F](y)=F|_{\mathrm{diag},+}(y)-F|_{\mathrm{diag},-}(y).
\]
Thus $+$ means approach from $x>y$, and $[F]$ denotes a jump. Operator commutators are written $[A,B]=AB-BA$.
The full diagonal quantities are
\[
\rho_z=G_z|_{\mathrm{diag}},\qquad
\sigma_z=(\partial_x-\partial_y)G_z|_{\mathrm{diag}},\qquad
\tau_z=(\partial_x-\partial_y)^2G_z|_{\mathrm{diag},+}.
\]
The first two are independent of the side; the second derivatives entering $\tau_z$ may jump. After subtracting the free terms, write
\[
r_z=\rho_z-\rho_{0,z},\qquad s_z=\sigma_z-\sigma_{0,z}.
\]
Here $\rho_{0,z},\sigma_{0,z},\tau_{0,z}$ are obtained from the preceding definitions by replacing $G_z(x,y)$ with $g_z(x-y)$; they are spatial constants.
At the reference point $\eta=z_1$, Lemma~\ref{lem:variation}
identifies these quantities with the variational coordinates:
\[
\mathcal G(U)=(r_\eta,s_\eta)
\]
for Schwartz data.
At rough regularity, only $r_z=4\sqrt3g_v$ and $s_z=4g_u$, defined by distributional variation, will be used.

The formal resolvent identities will be used on the kernel diagonal. We first justify the classical domains and the regularity of the free-subtracted kernel, including the jumps at the diagonal.

\begin{lemma}\label{lem:kernel}
Fix an admissible real $z$ and Schwartz data in the spectral ball. Write $R=R_z$, $\Delta R=R-R_{0,z}$ and $\Delta G(x,y)=G_z(x,y)-g_z(x-y)$. Then
\[
R:L^2\to H^3,\qquad R:H^j\to H^{j+3}\quad(j\ge0),
\qquad R:H^{-2}\to H^1,\qquad \Delta R\in\mathfrak S_1,
\]
\[
\langle(x,y)\rangle^N(\partial_x+\partial_y)^j\Delta G\in H^2(\mathbb R^2)
\quad(j,N\ge0),\qquad r_z,s_z\in\mathcal S(\mathbb R).
\]
Near the diagonal,
\[
\Delta G\in C^3_{\rm loc},\qquad [G_{xx},G_{xy},G_{yy}]=(1,-1,1),\qquad
[\partial_x^j\partial_y^{3-j}G]=0\quad(0\le j\le3).
\]
For every smooth finite-dimensional family of admissible spectral points and Schwartz data in the spectral ball, $\Delta R_z(U)$ is smooth in $\mathfrak S_1$, including mixed spectral and data derivatives.
The constants may depend on the fixed spectral point, inverse margin, and finitely many Schwartz seminorms.
\end{lemma}

\begin{proof}
We establish the inverse mapping properties and weighted kernel bounds, then the diagonal jumps and $C^3$ regularity of the difference kernel. The final factorization proves trace-class membership and trace-norm parameter dependence.

\emph{Matrix inverse.} On \(\mathcal H=L^2(\mathbb R;\mathbb C^3)\), define
\[
\mathcal A_z=\partial_x I_3-M_z,\qquad\operatorname{Dom}\mathcal A_z=H^1(\mathbb R;\mathbb C^3),
\qquad M_z=\begin{pmatrix}0&1&0\\0&0&1\\z-\mathfrak a'/2-\mathfrak b&-\mathfrak a&0\end{pmatrix}.
\]
The graph norm is equivalent to \(H^1\). With \(\mathcal E=M_z-M_{0,z}\),
\[
\mathcal E_{31}=-u'/4-v/(4\sqrt3),\qquad
\mathcal E_{32}=-u/2,\qquad\mathcal E_{ij}=0\text{ otherwise}.
\]
In this proof, write $M=M_z$, $\mathcal A=\mathcal A_z$, $M_0=M_{0,z}$, and $\mathcal A_0=\partial_x I_3-M_0$.
For the roots \(r_j\) of \(r^3-r/4=z\), put
\[
P_j=\prod_{\ell\ne j}\frac{M_0-r_\ell I}{r_j-r_\ell},\qquad
K_0(s)=\mathbf1_{s>0}\sum_{\Re r_j<0}e^{r_js}P_j
-\mathbf1_{s<0}\sum_{\Re r_j>0}e^{r_js}P_j.
\]
Then
\[
[K_0]=I,\qquad\mathcal A_0K_0=K_0\mathcal A_0=I,\qquad
|K_0(s)|\le C_ze^{-c_z|s|},\qquad K_0:L^2\to H^1.
\]
The weak scalar inverse satisfies, for \(f\in L^2\),
\[
Rf=R_0(f-VRf)\in H^3,\qquad (L-z)Rf=f,\qquad R(L-z)\psi=\psi\quad(\psi\in H^3).
\]
Thus the balanced inverse is the inverse of the classical realization with domain \(H^3\).

Split \(\mathcal E=\mathcal E_c+\mathcal E_t\) using a cutoff supported in \(I_c=[-2R,2R]\). For \(0<\delta<c_z\), equip kernels with
\[
\|K\|_\delta=\mathop{\mathrm{ess\,sup}}_{x,y}e^{\delta|x-y|}|K(x,y)|.
\]
Choose \(R\) so that the following contraction and its \(L^2\) analogue have norm at most \(1/2\):
\begin{equation}\label{eq:legacy-24}
\|K_0\mathcal E_tK\|_\delta
\le C_z\|\mathcal E_t\|_\infty
\int e^{-(c_z-\delta)|s|}\,ds\,\|K\|_\delta.
\end{equation}
Then
\[
K_t=\sum_{n\ge0}(K_0\mathcal E_t)^nK_0,\qquad
\|K_t\|_\delta\le2\|K_0\|_\delta,\qquad
\mathcal A_tK_t=K_t\mathcal A_t=I.
\]
The series converges both as kernels and in \(\mathcal B(L^2)\), since
\(\|K\|_{\mathrm{op}}\le2\delta^{-1}\|K\|_\delta\).

Set \(C=\mathcal E_cK_t\). Exponential decay gives
\[
\|C\|_{\mathfrak S_2}^2\lesssim_{z,\delta}\|\mathcal E_c\|_2^2.
\]
If \(F=CF\), then \(Y=K_tF\in H^1\) satisfies \(\mathcal A Y=0\). Its first two rows give
\[
Y=(\psi,\psi',\psi''),\qquad\psi\in H^3,\qquad(L-z)\psi=0.
\]
The scalar inverse implies \(F=0\). Consequently \(Q=(I-C)^{-1}\in\mathcal B(\mathcal H)\). Moreover,
\begin{equation}\label{eq:legacy-25}
Q(\mathbf1_{I_c}\mathcal H)\subset\mathbf1_{I_c}\mathcal H,\qquad
K=K_tQ=K_t+K_tQ\mathcal E_cK_t,\qquad\mathcal AK=K\mathcal A=I.
\end{equation}
The correction kernel is the row-column \(L^2\) pairing
\[
K_{\mathrm{corr}}(x,y)=\int_{I_c}K_t(x,s)
\,[Q(\mathcal E_c(\cdot)K_t(\cdot,y))](s)\,ds.
\]
Its factors obey
\[
\begin{aligned}
\|K_t(x,\cdot)\mathbf1_{I_c}\|_2&\lesssim |I_c|^{1/2}e^{-\delta d(x,I_c)},\\
\|\mathcal E_cK_t(\cdot,y)\|_2&\lesssim\|\mathcal E_c\|_2e^{-\delta d(y,I_c)}.
\end{aligned}
\]
Hence
\[
|K_{\mathrm{corr}}(x,y)|\lesssim\|Q\|e^{-\delta(d(x,I_c)+d(y,I_c))},\qquad
|K(x,y)|\lesssim e^{-\delta|x-y|}.
\]
Strong continuity of the rows and columns in \(L^2\) makes the correction continuous. The only jump of \(K\) is \(I\), and
\[
K_x=M(x)K+I\delta(x-y),\qquad K_y=-KM(y)-I\delta(x-y).
\]

\emph{Scalar kernel and weighted derivatives.} The operator with kernel \(G=e_1^TKe_3\) is the scalar inverse \(R\). For \(\Delta K=K-K_0\),
\[
\Delta K=K_0\mathcal EK,\qquad
\Delta K_x=M_0\Delta K+\mathcal EK,\qquad
\Delta K_y=-\Delta KM-K_0\mathcal E.
\]
The full second derivatives are
\begin{equation}\label{eq:legacy-26}
\begin{aligned}
\Delta K_{xx}&=M_0^2\Delta K+M_0\mathcal EK+\mathcal E'K+\mathcal EMK+\mathcal E(x)\delta,\\
\Delta K_{xy}&=-M_0\Delta KM-M_0K_0\mathcal E-\mathcal EKM-\mathcal E(x)\delta,\\
\Delta K_{yy}&=\Delta KM^2-\Delta KM'+K_0\mathcal EM+K_0M_0\mathcal E-K_0\mathcal E'+\mathcal E(y)\delta.
\end{aligned}
\end{equation}
The scalar projection cancels every displayed atom because \(\mathcal E_{13}=0\). For \(\delta'<\delta\),
\[
\langle x\rangle^m\langle y\rangle^m|\Delta K(x,y)|
\lesssim_m\int e^{-\delta'|x-s|}\langle s\rangle^{2m}|\mathcal E(s)|e^{-\delta'|s-y|}\,ds,
\]
whose \(L^2_{x,y}\) norm is bounded by \(C\|\langle s\rangle^{2m}\mathcal E\|_1\). Equation \eqref{eq:legacy-26} gives the same weighted \(L^2\) bounds for scalar derivatives through order two.

For \(\mathcal C=\partial_x+\partial_y\), differentiation of the translated inverse yields
\begin{equation}\label{eq:legacy-27}
\mathcal C^jK=
\sum_{\substack{n_1+\cdots+n_\ell=j\\n_i\ge1}}
\frac{j!}{n_1!\cdots n_\ell!}
K\mathcal E^{(n_1)}K\cdots\mathcal E^{(n_\ell)}K,\qquad j\ge1.
\end{equation}
Each term contains a Schwartz coefficient. Using a smaller exponential rate for each fixed finite convolution gives
\[
\langle x\rangle^m\langle y\rangle^m\mathcal C^j\Delta G\in H^2(\mathbb R^2)
\qquad(j,m\ge0).
\]
Here \(\mathcal C\delta(x-y)=0\), and the differentiated scalar atoms still vanish. The diagonal trace theorem and all center derivatives imply
\[
r=\Delta G|_{x=y}\in\mathcal S,\qquad
s=(\partial_x-\partial_y)\Delta G|_{x=y}\in\mathcal S.
\]

\emph{Local diagonal regularity.} Let \(F'=MF\) and
\[
P(x,y)=\mathbf1_{x>y}e_1^TF(x)F(y)^{-1}e_3.
\]
Zero initial values in the scalar ODE give both inverse identities on compact tests. The difference \(G-P\) is homogeneous in each variable; expanding in a fundamental basis in \(x\) leaves coefficients satisfying the smooth transposed ODE in \(y\). Thus \(G-P\in C^\infty\). For \(e=x-y>0\),
\[
e_1^TF(y+e)F(y)^{-1}e_3
=\frac{e^2}{2}-\frac{\mathfrak a(y)e^4}{24}+O(e^5).
\]
The free subtraction starts at \(\mathbf1_{e>0}e^4\). Therefore
\[
\Delta G\in C^3_{\mathrm{loc}},\qquad
[G_{xx},G_{xy},G_{yy}]=(1,-1,1),\qquad
[\partial_x^j\partial_y^{3-j}G]=0\quad(0\le j\le3).
\]
These side values, rather than second-order \(H^2\) traces, justify the normal diagonals below.

\emph{Trace-class resolvent difference and parameter dependence.} Fix \(w(x)=\langle x\rangle^{-2}\). For
\(f_1=u/2\), \(f_0=u'/4+\mathfrak b\), and \(C_U=I-VR\in\mathcal B(L^2)\),
\begin{equation}\label{eq:legacy-28}
\Delta R=-\sum_{j=0}^1(R_0M_w)(M_{f_j/w}\partial_x^jR_0)C_U.
\end{equation}
Both displayed kernel factors are Hilbert--Schmidt:
\[
\|R_0M_w\|_{\mathfrak S_2}=\|g_z\|_2\|w\|_2,\qquad
\|M_{f_j/w}\partial_x^jR_0\|_{\mathfrak S_2}=\|f_j/w\|_2\|g_z^{(j)}\|_2.
\]
Since $C_U$ is bounded, \eqref{eq:legacy-28} proves $\Delta R\in\mathfrak S_1$.
The inverse map in \(\mathcal B(L^2,H^3)\), the free multipliers, and \(f_j/w\in L^2\) are smooth on compact smooth parameter families. Hence \eqref{eq:legacy-28} is smooth in \(\mathfrak S_1\), including mixed spectral/data derivatives.

The same elliptic estimates for the transposed operator give $R:H^{-2}\to H^1$. This establishes all the asserted mapping, kernel, and trace-norm properties.
This completes the proof.
\end{proof}

The regularized determinant variation can now be identified with the subtracted Green diagonals. The associated boundary operators are needed for the two-parameter comparison.

\begin{lemma}\label{lem:variation}
For Schwartz data in the spectral ball and $h=(h_u,h_v)\in C_c^\infty(\mathbb R;\mathbb R^2)$,
\[
V_h=\tfrac14(M_{h_u}\partial_x+\partial_x\circ M_{h_u})+\tfrac1{4\sqrt3}M_{h_v},
\qquad dA_z(U)[h]=\operatorname{Tr}(\Delta R_zV_h),
\]
\[
\frac{\delta A_z}{\delta u}=\frac{s_z}{4},\qquad
\frac{\delta A_z}{\delta v}=\frac{r_z}{4\sqrt3}.
\]
The operator with kernel $\Pi_z(x,y)=\mathbf1_{x>y}G_z(x,y)$ preserves $H^3$. Its boundary identities on $H^3$, with values in $L^2$, are
\begin{equation}\label{eq:legacy-30}
\begin{gathered}
(L-z)\Pi_z=Q_z^\ell,\qquad \Pi_z(L-z)=Q_z^r,\\
Q_z^\ell=\rho_z\partial_x^2+\tfrac12(5\rho_z'+\sigma_z)\partial_x 
+\tfrac14(7\rho_z''+4\sigma_z'+\tau_z)+\mathfrak a\rho_z,\\
Q_z^r=\rho_z\partial_x^2+\tfrac12(-\rho_z'+\sigma_z)\partial_x 
+\tfrac14(\rho_z''-2\sigma_z'+\tau_z)+\mathfrak a\rho_z.
\end{gathered}
\end{equation}
\end{lemma}

The second identity compares two spectral parameters. Its diagonal evaluations give the complete transformed auxiliary field, retaining the free terms as well as the perturbation.

\begin{lemma}\label{lem:two-parameter}
Let $U$ be Schwartz data in the spectral ball and let $\zeta\ne\eta$ be admissible real spectral points. Then
\begin{equation}\label{eq:legacy-31}
(\zeta-\eta)[\Pi_\zeta,R_\eta]=R_\eta Q_\zeta^\ell-Q_\zeta^rR_\eta,
\end{equation}
as an identity in $\mathcal L(L^2,H^1)$. Put
\[
(\rho,\sigma,\tau)=(\rho_\eta,\sigma_\eta,\tau_\eta),\qquad
(c,d,f)=(\rho_\zeta,\sigma_\zeta,\tau_\zeta),
\]
and let $K$ be the kernel of $R_\eta Q_\zeta^\ell-Q_\zeta^rR_\eta$. Its zeroth and first normal diagonals satisfy
\begin{equation}\label{eq:legacy-33}
K|_{\mathrm{diag}}=-c\sigma'+c'\sigma/2-d\rho'/2+d'\rho.
\end{equation}
\begin{equation}\label{eq:legacy-34}
(\partial_x-\partial_y)K|_{\mathrm{diag}}=-c\tau'+f'\rho
+(c''\rho'-c'\rho'')/2+(d'\sigma-d\sigma')/2.
\end{equation}
For the normalized field, fix $\eta=z_1$, let $k\ge2$, and set
\[
\begin{gathered}
H_k=C_ks_k,\qquad Y_k=J\nabla H_k,\\
p_k=\frac{C_k}{48\sqrt3}\sum_{\zeta=\pm z_k}\frac{\rho_\zeta}{\zeta-\eta},\qquad
q_k=\frac{C_k}{48\sqrt3}\sum_{\zeta=\pm z_k}\frac{\sigma_\zeta}{\zeta-\eta},\\
t_k=-\tfrac13[p_k'''+4(\mathfrak a p_k)'],\qquad
t_\eta=-\tfrac13[\rho'''+4(\mathfrak a\rho)']=\tau'.
\end{gathered}
\]
Then the Green map \eqref{eq:reference-green-map} satisfies $D_U \mathcal G(U)Y_k(U)=(F_{k,1},F_{k,2})$, where
\begin{equation}\label{eq:legacy-35}
\begin{aligned}
F_{k,1}&=-p_k\sigma'+\tfrac12p_k'\sigma-\tfrac12q_k\rho'+q_k'\rho,\\
F_{k,2}&=-p_kt_\eta+t_k\rho
+\tfrac12(p_k''\rho'-p_k'\rho'')
+\tfrac12(q_k'\sigma-q_k\sigma').
\end{aligned}
\end{equation}
\end{lemma}

\begin{proof}[Proof of Lemmas~\ref{lem:variation} and~\ref{lem:two-parameter}]
We first prove the variational identities and the two boundary formulas \eqref{eq:legacy-30}. We then derive the commutator identity \eqref{eq:legacy-31}, its diagonal evaluations \eqref{eq:legacy-33}--\eqref{eq:legacy-34}, and the normalized field \eqref{eq:legacy-35}.

For compact smooth $h$, the localized factorization of $\Delta R\,V_h$ uses kernel derivatives through order two. The weighted bounds in Lemma~\ref{lem:kernel} give Hilbert--Schmidt factors, so $\Delta R\,V_h\in\mathfrak S_1$ and its trace is an absolutely convergent product-kernel integral. The abstract determinant variation in Lemma~\ref{lem:determinant}, followed by trace cyclicity, therefore gives
\[
dA_z(U)[h]=\operatorname{Tr}(\Delta R\,V_h).
\]
The trace-norm parameter dependence from Lemma~\ref{lem:kernel} also permits mixed spectral and data differentiation. The local $C^3$ regularity of $K=\Delta G$ justifies diagonal integration by parts:
\[
\operatorname{Tr}(K M_h\partial_x)=-\int(h'K+hK_y)|_{\rm diag},\qquad
\operatorname{Tr}(K\partial_x\circ M_h)=-\int hK_y|_{\rm diag}.
\]
Integrating the $h'$ term once more, the two $u$ vertices give \(\frac14\int h(K_x-K_y)|_{\rm diag}\) and the $v$ vertex gives \(\frac1{4\sqrt3}\int hK|_{\rm diag}\). Thus the gradient is \(g_u=s_z/4,g_v=r_z/(4\sqrt3)\). Rough coordinates remain defined by distributional variation.

By Lemma~\ref{lem:kernel},
\begin{equation}\label{eq:legacy-29}
[G_{xx},G_{xy},G_{yy}]=(1,-1,1),\qquad
[\partial_x^j\partial_y^{3-j}G]=0\quad(0\le j\le3),
\end{equation}
and \(G-g_z\in C^3_{\mathrm{loc}}\). All normal diagonals below use the positive-side extension.

The kernel of $\Pi_z$ is $\mathbf1_{x>y}G_z(x,y)$, with the positive-side convention fixed above.

The positive-side homogeneous equation gives the distributional boundary:\[
(L-z)\Pi_z:\quad
\delta''G+3\delta'G_x+\delta(3G_{xx}+\mathfrak aG).
\]
Using \(f\delta'=f_d\delta'+(f_y)_d\delta\) and \(f\delta''=f_d\delta''+2(f_y)_d\delta'+(f_{yy})_d\delta\),
the coefficients are \(\rho,\ 3G_x+2G_y,\ 3G_{xx}+3G_{xy}+G_{yy}+\mathfrak a\rho\).
Integration by parts against compact tests on the right gives the coefficients \(\rho,-G_y,G_{yy}+\mathfrak a\rho\).
Substituting\[\begin{gathered}
G_x=(\rho'+\sigma)/2,\\ G_y=(\rho'-\sigma)/2,\\
G_{xx}=(\rho''+2\sigma'+\tau)/4,\\
G_{xy}=(\rho''-\tau)/4,\\
G_{yy}=(\rho''-2\sigma'+\tau)/4
\end{gathered}\]
gives precisely the coefficients in \eqref{eq:legacy-30}.
We always use the positive-side extension and do not evaluate \(\theta \delta \).
After the test identity is established, \(\Pi \) is a free \(L^1\) convolution plus a Hilbert--Schmidt difference kernel. Thus $\Pi_z$ is bounded on $L^2$. Equation \eqref{eq:legacy-30} and elliptic estimates show that \(\Pi \) preserves \(H^3\). Since \(Q^\ell:L^2\to H^{-2}, R:H^{-2}\to H^1,Q^rR:L^2\to H^1\), products are defined on a dense subspace of \(H^3\).
For distinct \(\zeta\ne\eta\), we obtain \eqref{eq:legacy-31}.
The identity is first proved on tests. The bounds $Q_\zeta^\ell:L^2\to H^{-2}$, $R_\eta:H^{-2}\to H^1$, and $Q_\zeta^rR_\eta:L^2\to H^1$ extend it to the asserted identity in $\mathcal L(L^2,H^1)$. No rough diagonal is taken at this stage.

Let the reference full diagonal be \((\rho,\sigma,\tau)\) and the other one \((c,d,f)\). Here $d$ is a diagonal function, not the spectral polynomial \(d_k\).
The kernel on the right of \eqref{eq:legacy-31} is\[
K(x,y)=c(y)G_{yy}-c(x)G_{xx}
+a_1(y)G_y-a_2(x)G_x+(a_3(y)-a_4(x))G,
\]\begin{equation}\label{eq:legacy-32}
a_1=(-c'-d)/2,\ a_2=(-c'+d)/2,\quad
a_3=(c''+2d'+f)/4+\mathfrak a c,\quad a_4=(c''-2d'+f)/4+\mathfrak a c.
\end{equation}
The jumps cancel:
\[
[K]=c(1-1)=0,\qquad
[K_x]=-c'-a_1-a_2=0,\qquad[K_y]=c'+a_1+a_2=0.
\]
Equation \eqref{eq:legacy-29} gives \(K\in C^1\). For the left commutator, the boundary is
\[
(\rho_\zeta(x)-\rho_\zeta(y))G_\eta(x,y),
\]
which vanishes at \(x=y\); row-column \(L^2\) bounds justify differentiation. Equality of the distribution kernels therefore implies equality of their zeroth and first normal diagonals.

Evaluation at the diagonal gives \eqref{eq:legacy-33}.
Before simplification, the second diagonal is\[
\begin{split}
(\partial_x-\partial_y)K|_{\mathrm{diag}}={}&-c\tau'-c'(\rho''+\tau)/2\\
&+(a_1-a_2)\sigma'/2-(a_1+a_2)\tau/2\\
&-(a_1'+a_2')\rho'/2+(a_1'-a_2')\sigma/2\\
&+d'\sigma-(a_3'+a_4')\rho .
\end{split}
\]
The left and right homogeneous equations give \(c'''+3f'=-4(\mathfrak a c)'\), hence \(a_3'+a_4'=-f'\).
Combining the terms yields \eqref{eq:legacy-34}.

The boundary formulas \eqref{eq:legacy-30} also give $L_{J\nabla A_z}=[-\Pi_z/(24\sqrt3),L]$.
After normalization, \(R_{\eta,Y_k}=[C_k(\Pi_++\Pi_-)/(48\sqrt3),R_\eta]\).
Summing \eqref{eq:legacy-33}--\eqref{eq:legacy-34} with the coefficients defining $p_k,q_k,t_k$ gives the normalized field \eqref{eq:legacy-35}.
The variational formula and boundary identities prove Lemma~\ref{lem:variation}. Equations \eqref{eq:legacy-31} and \eqref{eq:legacy-33}--\eqref{eq:legacy-35} prove Lemma~\ref{lem:two-parameter}. This completes both proofs.
\end{proof}

\begin{remark}\label{rem:free-field-coefficients}
Write $p_k^0=p_k(0)$ and $q_k^0=q_k(0)$. The free diagonal values $\rho_{0,z_k}=-1/d_k$, $\rho_{0,-z_k}=1/d_k$, and $\sigma_{0,\pm z_k}=-2k/d_k$ give\begin{equation}\label{eq:legacy-36}
p_k^0=-\frac{C_kz_k}{24\sqrt3d_k(z_k^2-\eta^2)},\qquad
q_k^0=-\frac{C_kk\eta}{12\sqrt3d_k(z_k^2-\eta^2)}.
\end{equation}
For \(k\ge 2\), the bounds are respectively $O(1)$ and \(O(k^{-2})\), uniformly. The free $p$ term converges to \(-\sqrt{3}\). The free terms must be kept, and the two \(\sigma \) terms are estimated together.
\end{remark}

\subsection{Poisson commutation at the endpoint}\label{subsec:poisson-commutation-at-the-endpoint}

The two-parameter identities imply conservation of one spectral functional under the Hamiltonian field of another. Continuity of the abstract pairing then extends this conclusion to the endpoint ball.

\begin{lemma}\label{lem:poisson}
For $U$ in the real endpoint spectral ball and all admissible real $z,\zeta$,
\[
\{A_\zeta,A_z\}(U)
=\langle\nabla A_\zeta(U),J\nabla A_z(U)\rangle=0.
\]
The pairing is the continuous $X^*$--$X$ pairing, including equal and opposite spectral points. The corresponding local Hamiltonian flows commute whenever defined.
\end{lemma}

\begin{proof}
We first show that $dA_\zeta(U)[J\nabla A_z(U)]=0$ for Schwartz data. We then extend this identity to the endpoint ball and deduce commutation of the local Hamiltonian flows.

For Schwartz data, Lemma~\ref{lem:variation} gives
\[
Y_z=J\nabla A_z,\qquad
\delta\mathfrak a=r_z'/(8\sqrt3),\quad\delta\mathfrak b=s_z'/(16\sqrt3),\qquad
L_{Y_z}=[\mathcal T,L],\quad\mathcal T=-\Pi_z/(24\sqrt3).
\]
On the proper inverse domains,
\[
R_{\zeta,Y_z}=[\mathcal T,R_\zeta].
\]
Write \(\Pi=\Pi_0+\Delta\Pi\), \(R=R_0+\Delta R\). Free convolutions commute; the bounded/trace-class commutators have zero trace. For the cross terms,
\[
\begin{aligned}
\Delta\Pi R_0&=(\Delta\Pi\langle x\rangle)(\langle x\rangle^{-1}R_0),\\
R_0\Delta\Pi&=(R_0\langle x\rangle^{-1})(\langle x\rangle\Delta\Pi).
\end{aligned}
\]
Each factor is Hilbert--Schmidt by Lemma~\ref{lem:kernel}, and both traces equal the same absolutely convergent double kernel integral. Hence
\[
\operatorname{Tr}[\mathcal T,R_\zeta]=0.
\]

Spectral differentiation gives \(\partial _\zeta T=R_{0}T/2\) and \(\partial _\zeta B=(R_{0}B+BR_{0})/2\), hence\begin{equation}\label{eq:legacy-37}
\partial_\zeta A_\zeta=-\operatorname{Tr}(R_0S B^2)
=\operatorname{Tr}(R_0-R_\zeta)-\operatorname{Tr}(R_0V R_0).
\end{equation}
This identity follows from
\[
\Delta R=-TBT+TB^2ST
\]
in trace-class differences. Mixed differentiation is permitted by Lemma~\ref{lem:kernel}. For \(b_z(\zeta)=dA_\zeta[Y_z]\),
\[
b_z'(\zeta)=-\operatorname{Tr}[\mathcal T,R_\zeta]
-\operatorname{Tr}(R_0V_{Y_z}R_0)=0.
\]
The last free trace is a linear combination of
\[
\int r_z'\,dx=\int r_z''\,dx=\int s_z'\,dx=0.
\]
Thus \(b_z\) is constant on each real ray. Estimate \eqref{eq:legacy-8} gives
\[
|b_z(\zeta)|\lesssim |\zeta|^{-2/3}\|U\|_X\|Y_z\|_X\to0
\quad(|\zeta|\to\infty),
\]
so the bracket vanishes, including equal and opposite spectral points. Density and continuity extend it to \(X\). Symmetric Hessians and the constant skew operator \(J\) imply
\[
D_U Y_\zeta Y_z-D_U Y_zY_\zeta=0.
\]
Uniqueness in the Banach ODE variation equation proves local commutation.
This proves the Poisson identity and the stated local commutation.
\end{proof}

\section{Auxiliary flows and commutation}\label{sec:flows}
The conserved spectral bound gives global auxiliary flows and global high-regularity original flows. We prove their commutation in a common variation class, with the generator domains stated explicitly. This supplies the chain rule needed for the difference integral.

\subsection{The auxiliary Hamiltonian flows}\label{subsec:the-auxiliary-hamiltonian-flows}

The spectral average defines the regularized Hamiltonians on the endpoint space. Their fixed-parameter fields are ordinary Banach-space vector fields; the common spectral barrier yields their global auxiliary times.

\begin{lemma}\label{lem:auxiliary}
Let $H_k=C_ks_k$ and $k\ge1$. On the endpoint spectral ball,
\[
Y_k=J\nabla H_k:X\to X^{1/2}\text{ is analytic},\qquad
\|Y_k(U)\|_{X^{1/2}}+\|D_U Y_k(U)\|_{\mathcal L(X,X^{1/2})}
\le C_{\beta_*}k^3.
\]
Its flow $\Psi_k^a$ is defined on $\mathcal D$ for every $a\in\mathbb R$ and
\[
A_z(\Psi_k^aU)=A_z(U),\qquad
\mathscr S_R(\Psi_k^aU)=\mathscr S_R(U).
\]
The auxiliary flows commute. Higher-regularity and finite-time data estimates at fixed $k$ may depend on $k$.
\end{lemma}

\begin{proof}
We prove the fixed-$k$ vector-field bound and spectral conservation, use the spectral barrier to obtain global auxiliary time, and finally extend local commutation along the global trajectories.

Lemmas~\ref{lem:determinant} and~\ref{lem:poisson} give
\[
\|Y_k(U)\|_{X^{1/2}}+\|D_U Y_k(U)\|_{\mathcal L(X,X^{1/2})}
\le C_{\beta_*}k^3,\qquad
\partial_tA_\zeta(\Psi_k^tU)=dA_\zeta[Y_k]=0.
\]
Thus \(\mathscr S\) and \(\mathscr S_R\) are conserved, and
\[
\gamma\|\Psi_k^tU\|_X^2\le\mathscr S(U)<\gamma\varepsilon^2/4.
\]
The field is bounded and Lipschitz on \(\overline B_X(0,3\varepsilon/4)\). A finite maximal time would give a strong \(X\) limit in \(\mathcal D\), where the Banach ODE restarts. Local commutation extends along the global trajectories.
The field bound, global auxiliary flow, conservation laws, and commutation are therefore established. This completes the proof.
\end{proof}

\subsection{The high-regularity original flow and variations}\label{subsec:the-high-regularity-original-flow-and-variations}

For comparison with the auxiliary dynamics we need an original flow on smooth data, together with its high-space variations. Spectral conservation and the low-norm lifespan give this flow before any rough original limit is taken.

\begin{lemma}\label{lem:high-flow}
For $U_0\in\mathcal D\cap X^{7/2}$, the original flow is defined for all real times and
\[
\Phi_E^\cdot U_0\in C(\mathbb R;X^{7/2})\cap C^1(\mathbb R;X^{3/2}),
\qquad \Phi_E^tU_0\in\mathcal D.
\]
For $m=3/2,7/2,11/2$, the original local solution map is $C^2$ on $X^m$ over its interval of existence. Writing $U(t)=\Phi_E^tU_0$, its first and second variations $Z_h=D_U \Phi_E^\cdot(U_0)h$ and $K_{h,l}=D_U^2 \Phi_E^\cdot(U_0)[h,l]$ are the unique $C_tX^m$ solutions of
\begin{equation}\label{eq:legacy-41}
Z_h(t)=\mathcal U_E(t)h+\int_0^t\mathcal U_E(t-a)D_U N(U(a))Z_h(a)\,da.
\end{equation}
\begin{equation}\label{eq:legacy-42}
K_{h,l}(t)=\int_0^t\mathcal U_E(t-a)
\{D_U N(U(a))K_{h,l}(a)+D_U^2 N[Z_h(a),Z_l(a)]\}\,da.
\end{equation} At fixed $k$, the auxiliary maps satisfy
\[
Y_k:X^{n-1/2}\longrightarrow X^{n+1/2},\qquad
\Psi_k^a:X^{n-1/2}\longrightarrow X^{n-1/2},\qquad n=0,\ldots,6,
\]
with $C^2$ dependence on the data on their spectral domains. Along an auxiliary trajectory $U(a)=\Psi_k^aU_0$, the first and second variations $P_h,Q_{h,l}$ are the unique solutions on each finite interval of
\begin{equation}\label{eq:legacy-45}
P_h'=D_U Y_k(U)P_h,\quad P_h(0)=h;
\qquad Q_{h,l}'=D_U Y_k(U)Q_{h,l}+D_U^2 Y_k(U)[P_h,P_l],\quad Q(0)=0.
\end{equation} Original data in $\mathcal D\cap X^{11/2}$ retain that regularity on every finite time interval.
\end{lemma}

\begin{proof}
We first prove global continuation on $\mathcal D\cap X^{7/2}$. We then establish the generator domains, the first and second variation equations, and the asserted high-regularity properties of both flows.

Before a first exit from \(B_X(0,3\varepsilon/4)\), conservation gives
\[
\|\Phi_E^tU_0\|_X^2\le\gamma^{-1}\mathscr S(U_0)<\varepsilon^2/4.
\]
Lemma~\ref{lem:persistence} supplies the fixed lifespan \(\tau_0\). Restarting at a time less than \(\tau_0/2\) from a putative finite endpoint extends the solution. High mild uniqueness joins the intervals. No terminal \(X\) limit or uniform high norm is required.

\smallskip\noindent\emph{Spaces and domains.} For the auxiliary regularity argument, put
\[
E_n=X^{n-1/2},\qquad
\mathcal Q_n^\nabla=H^{n+1/2}\times H^{n+3/2}.
\]
Recall the decomposition $\mathcal V_E(U)=\mathscr L_EU+N(U)$ from Section~\ref{subsec:notation}, with $\mathscr L_E$ defined in \eqref{eq:linear-flow-notation}. The two-derivative loss is recorded by
\begin{equation}\label{eq:generator-domains}
\operatorname{Dom}(\mathscr L_E\text{ on }X^m)=X^{m+2},\qquad
N\in C^\infty(X^m;X^m),\quad m>1/2.
\end{equation}
For $U_0\in X^{7/2}$, the original trajectory and its variations satisfy
\[
\begin{aligned}
U&\in C_tX^{7/2}\cap C_t^1X^{3/2},\\
h\in X^{3/2}&\Longrightarrow Z_h\in C_tX^{3/2}\cap C_t^1X,\\
h\in X^{7/2}&\Longrightarrow Z_h\in C_tX^{7/2}\cap C_t^1X^{3/2}.
\end{aligned}
\]
At fixed $k$, the gradient and auxiliary-field estimates have the form
\[
\nabla H_k:E_n\longrightarrow\mathcal Q_n^\nabla,\qquad
Y_k:E_n\longrightarrow E_{n+1},
\]
with the first and second directional derivatives in the same corresponding output spaces. In particular,
\[
Y_k(U(t))\in C_tX^{9/2}\cap C_t^1X^{5/2}.
\]
High variations agree with base variations by Duhamel uniqueness. Constants in this high-regularity argument may depend on the fixed \(k\).

\smallskip\noindent\emph{The linear group and nonlinearity.} Set $\varpi(\xi)=\xi\langle\xi\rangle$. Under the isometry $\widehat U\mapsto(\widehat u,\widehat v/\langle\xi\rangle)$, the generator is $i\varpi\begin{pmatrix}0&1\\1&0\end{pmatrix}$. In the original variables the group is
\begin{equation}\label{eq:legacy-38}
\widehat{\mathcal U_E(t)U}(\xi)
=\begin{pmatrix}
\cos(t\varpi(\xi))&i\langle\xi\rangle^{-1}\sin(t\varpi(\xi))\\
i\langle\xi\rangle\sin(t\varpi(\xi))&\cos(t\varpi(\xi))
\end{pmatrix}\widehat U(\xi).
\end{equation}
Unitarity in the weighted variables and dominated convergence give strong continuity. Moreover,
\[
\operatorname{Dom}\mathscr L_E
=\{U:\xi\langle\xi\rangle\langle\xi\rangle^m
(\widehat u,\widehat v/\langle\xi\rangle)\in L^2\}=X^{m+2}.
\]
Reality follows from oddness of \(\xi\langle\xi\rangle\).

For $m>1/2$, weighted convolution and \(\|\widehat{f}\|_{1}\le C_m\|f\|_{H^m}\) give \(\|fg\|_{H^m}\le C_m\|f\|_{H^m}\|g\|_{H^m}\). Thus $N$ is a continuous quadratic map on \(X^m\to X^m\), with\begin{equation}\label{eq:legacy-39}
D_U N(U)[h]=(0,-2\partial_x (uh_u)),\qquad
D_U^2 N(U)[h,l]=(0,-2\partial_x (h_ul_u)),\quad D_U^3 N=0,
\end{equation}
and \(\|N(U)\|_{X^m}\le C_m\|U\|_{X^m}^2\), \(\|D_U N(U)\|\le 2C_m\|U\|_{X^m}\), \(\|D_U^2 N\|\le 2C_m\). The derivative loss of $\partial_x$ on \(H^m\to H^{m-1}\) is absorbed by the second component. These estimates apply at the high regularities used here, but not at $m=-1/2$.

\smallskip\noindent\emph{Data variations of the original flow.} For $m>1/2$, the mild equation is
\begin{equation}\label{eq:legacy-40}
U(t)=\mathcal U_E(t)U_0+\int_0^t\mathcal U_E(t-a)N(U(a))\,da.
\end{equation}
The algebra bound and \eqref{eq:legacy-39} give a $C^2$ local solution map by Picard iteration. On a fixed interval $[-T,T]$ of existence, put $M=1+\|U\|_{C_tX^m}$. The linearized Volterra equations \eqref{eq:legacy-41}--\eqref{eq:legacy-42} have unique $C_tX^m$ solutions and
\[
\|Z_h\|_{C_tX^m}\le e^{C_mMT}\|h\|_{X^m},\qquad
\|K_{h,l}\|_{C_tX^m}\lesssim_{m,M,T}\|h\|_{X^m}\|l\|_{X^m}.
\]
For $\Delta_h=\Phi_E^\cdot(U_0+h)-U$, the quadratic identity
\[
N(U+\Delta_h)-N(U)-D_U N(U)\Delta_h=(0,-\partial_x (\Delta_{h,u}^2))
\]
and Gronwall give
\[
\|\Delta_h-Z_h\|_{C_tX^m}\lesssim_{m,M,T}\|h\|_{X^m}^2.
\]
Differentiating the linearized equation once more gives $K_{h,l}$. Subtracting the corresponding equations proves continuity of both derivatives. The mild equations imply
\[
U,Z_h\in C_tX^m\cap C_t^1X^{m-2}.
\]
For higher data and directions, uniqueness of \eqref{eq:legacy-41} identifies the variations obtained in different Sobolev spaces. Conversely, differentiating $\mathcal U_E(-t)U(t)$ in $X^{m-2}$ recovers \eqref{eq:legacy-40}. All statements hold in both time directions. These variation bounds may depend on $M$; the lifespan used for continuation is supplied by Lemma~\ref{lem:persistence}.

\smallskip\noindent\emph{Higher regularity from translation equivariance.} Let \(F:E_0\to F_0\) be \(C^{n+3}\) and translation equivariant. For \(\partial_x^jU\in E_0\), \(0\le j\le n\), differentiating \(F(\tau_yU)=\tau_yF(U)\) gives
\begin{equation}\label{eq:legacy-43}
\partial_x^jF(U)=
\sum_{p=1}^j\ \sum_{\substack{j_1+\cdots+j_p=j\\j_i\ge1}}
\frac{j!}{p!j_1!\cdots j_p!}
D_U^pF(U)[\partial_x^{j_1}U,\ldots,\partial_x^{j_p}U].
\end{equation}
Distributional difference quotients identify these derivatives. Applying \eqref{eq:legacy-43} to
\[
(U,h)\mapsto D_U F(U)h,\qquad(U,h,\ell)\mapsto D_U^2 F(U)[h,\ell]
\]
and using the integral Taylor remainder gives the lifted \(C^2\) maps.

For \(F_{0}=X^{1/2}\), the sum of its spatial derivative norms through order $n$ is equivalent to \(X^{n+1/2}\); for \(F_{0}=\mathcal Q_0^\nabla \) it is equivalent to \(\mathcal Q_n^\nabla \). The full trace variation bounds in Lemma~\ref{lem:determinant} give \(C^\infty \) maps \(\nabla H_k\) on \(X\to \mathcal Q_0^\nabla \) and \(Y_k\) on \(X\to X^{1/2}\), with ball bounds \(M_{k,p}<\infty \) at each order. Translation conjugates $B$ unitarily and preserves traces, so \(H_k\) is translation invariant. Its gradient and $J$ are equivariant. Hence\begin{equation}\label{eq:legacy-44}
Y_k:E_n\to E_{n+1},\quad
\nabla H_k:E_n\to \mathcal Q_n^\nabla
\quad(n=0,\ldots,6)
\end{equation}
are \(C^2\), for directions in \(E_n\). The principle holds for fixed integer $n$; here only \(n\le 6\) is needed. Rough directions use $n=0$ throughout.

\smallskip\noindent\emph{High persistence on the low-norm interval.} Lemma~\ref{lem:persistence} and \eqref{eq:legacy-43} give
\[
U_0\mapsto U\in C_tX^m\cap C_t^1X^{m-2},\qquad m=3/2,7/2,11/2,
\]
with \(C^2\) data dependence on the same interval \(\tau_0\). The highest derivative is absorbed by the fixed local inverse; the lower terms are controlled inductively. For \(m=7/2,11/2\), Lemma~\ref{lem:conservation} applies. Thus
\[
\|U(t)\|_X<\varepsilon/2,\qquad
\tau_0=\lambda_*^{-2}
\]
at every restart. Finitely many restarts cover a bounded time window. The high variation bounds may depend on the high trajectory norm; \(\tau_0\) does not.

\smallskip\noindent\emph{Auxiliary high-regularity flows.} By Lemmas~\ref{lem:poisson} and~\ref{lem:auxiliary}, the base flow \(\Psi _k\) exists for all auxiliary times on \(\mathcal D \) and conserves all spectral functionals. Translations preserve \(\mathcal D \), and ODE uniqueness makes \(\Psi _k\) equivariant.

Fix a finite auxiliary interval \(|b|\le B_{0}\). On the larger Neumann ball, all derivatives of \(Y_k\) have bounds \(M_{k,p}\). Its base variations satisfy \eqref{eq:legacy-45}.
Integral contraction and Gronwall give
\[
\|P_h\|_X\le e^{M_{k,1}B_0}\|h\|_X,\qquad
\|Q_{h,\ell}\|_X\le M_{k,2}B_0e^{3M_{k,1}B_0}\|h\|_X\|\ell\|_X.
\]
Higher variations and translation equivariance, followed by \eqref{eq:legacy-43}, give
\[
\Psi_k^bU_0\in C_bE_n\cap C_b^1E_n,\qquad
\Psi_k^b:E_n\to E_n\text{ is }C^2.
\]
The highest spatial derivative solves an equation with linear term
\[
D_U Y_k(U)\partial_x^nU;
\]
all other terms involve lower derivatives. Inductive Gronwall gives persistence on every fixed \(k,B_0\) window. Uniqueness in \eqref{eq:legacy-45} identifies its lifted and base variations.
This proves global continuation, high persistence, and the stated first and second variation properties. This completes the proof.
\end{proof}

\subsection{Infinitesimal and actual commutation}\label{subsec:infinitesimal-and-actual-commutation}

A vanishing Poisson bracket alone is insufficient for the comparison of actual trajectories. We place the transported fields in one linearized Duhamel class and use its uniqueness to prove the next statement.

\begin{proposition}\label{prop:commutation}
For $U_0\in\mathcal D\cap X^{7/2}$ and fixed $k\ge1$,
\[
\Phi_E^t\Psi_k^aU_0=\Psi_k^a\Phi_E^tU_0,\qquad t,a\in\mathbb R.
\]
\end{proposition}

\begin{proof}
We prove the infinitesimal identity \eqref{eq:legacy-48}, obtain the transport identity \eqref{eq:legacy-49} by uniqueness of the linearized equation, and then deduce actual commutation \eqref{eq:legacy-50}.

On \(U\in \mathcal D \cap X^{7/2}\), the energy $E$ is \(C^3\) and\begin{equation}\label{eq:legacy-46}
g_E=(u-\partial_x^2u-u^2,v),\qquad
\mathcal H_E h=((1-\partial_x^2-2u)h_u,h_v).
\end{equation}
By \eqref{eq:legacy-44},
\[
\begin{gathered}
g_E\in H^{3/2}\times H^{5/2},\qquad\mathcal V_E\in X^{3/2},\\
g_H=\nabla H_k\in H^{9/2}\times H^{11/2},\qquad Y_k\in X^{9/2}.
\end{gathered}
\]
Smooth spectral conservation and density in \(X^{7/2}\) imply
\[
dH_k(U)[\mathcal V_E(U)]=0\qquad(U\in\mathcal D\cap X^{7/2}).
\]
Differentiate in a compact smooth direction \(h\). Hessian symmetry and skew symmetry of \(\partial_x \) give
\begin{equation}\label{eq:legacy-47}
0=d^2H_k[h,\mathcal V_E]+\langle g_H,J\mathcal H_Eh\rangle
=\langle D_U g_H[\mathcal V_E]-\mathcal H_EY_k,h\rangle.
\end{equation}
Both terms belong to \(H^{5/2}\times H^{7/2}\): the variation uses \eqref{eq:legacy-44} on \(X^{3/2}\), and \(\mathcal H_E\) loses two derivatives from \(Y_k\in X^{9/2}\). Therefore
\begin{equation}\label{eq:legacy-48}
D_U Y_k(U)\mathcal V_E(U)
=JD_U g_H[\mathcal V_E]
=J\mathcal H_EY_k
=D_U \mathcal V_E(U)Y_k(U)
\quad\text{in }X^{5/2}.
\end{equation}

For \(U(t)=\Phi_E^tU_0\), \(U_0\in\mathcal D\cap X^{7/2}\), put \(W(t)=Y_k(U(t))\). Then
\[
U\in C_tX^{7/2}\cap C_t^1X^{3/2},\qquad
W\in C_tX^{9/2}\cap C_t^1X^{5/2},\qquad
W'=\mathscr L_EW+D_U N(U)W.
\]
Since \(W\) takes values in the generator domain on \(X^{3/2}\),
\[
W(t)=\mathcal U_E(t)Y_k(U_0)
+\int_0^t\mathcal U_E(t-a)D_U N(U(a))W(a)\,da.
\]
Uniqueness in \eqref{eq:legacy-41} gives
\begin{equation}\label{eq:legacy-49}
Y_k(\Phi_E^tU_0)=D_U \Phi_E^t(U_0)Y_k(U_0)\quad\text{in }X^{3/2}.
\end{equation}
For \(C_t(b)=\Phi_E^t(\Psi_k^bU_0)\), the data chain rule gives
\[
\partial_bC_t(b)=D_U \Phi_E^t(\Psi_k^bU_0)Y_k(\Psi_k^bU_0)=Y_k(C_t(b)).
\]
Banach ODE uniqueness implies
\begin{equation}\label{eq:legacy-50}
\Phi_E^t\Psi_k^bU_0=\Psi_k^b\Phi_E^tU_0.
\end{equation}
The high data maps are continuous, and both flows preserve \(\mathcal D\). On any finite common rectangle their high norms are bounded by compactness of the fixed-parameter data curve. Local commutation thus extends over the rectangle.
Since the rectangle is arbitrary, \eqref{eq:legacy-50} holds for all real $t$ and $a$. This proves the proposition.
\end{proof}

\subsection{Common trajectories and the strong chain rule}\label{subsec:common-trajectories-and-the-strong-chain-rule}

Actual commutation removes the transported auxiliary field from the derivative of a common trajectory. The extra two derivatives in the data make that chain rule valid in the strong Green space.

\begin{lemma}\label{lem:chain-rule}
Let $U_0\in\mathcal D\cap X^{11/2}$, fix $k\ge1$ and $t\in\mathbb R$, and set
\[
V(a)=\Phi_E^a\Psi_k^{t-a}U_0,\qquad a\text{ between }0\text{ and }t.
\]
Then
\[
V\in C_aX^{11/2}\cap C_a^1X^{7/2},
\]
\begin{equation}\label{eq:legacy-52}
V'=\mathcal V_E(V)-Y_k(V)\quad\text{in }X^{7/2},
\end{equation}
\begin{equation}\label{eq:legacy-53}
\partial_a\mathcal G(V)
=D_U \mathcal G(V)[\mathcal V_E(V)-Y_k(V)]\quad\text{in }Z_G.
\end{equation}
Here $\mathcal G$ is the paired map \eqref{eq:reference-green-map}; for these data,
\[
\mathcal G(U)=(r_\eta(U),s_\eta(U))
=\bigl(4\sqrt3\,g_v(\eta;U),4g_u(\eta;U)\bigr).
\]
The integrated identity in $Z_G$ is
\begin{equation}\label{eq:legacy-54}
\begin{aligned}
&\mathcal G(\Phi_E^tU_0)-\mathcal G(\Psi_k^tU_0)\\
&\qquad=\int_0^tD_U \mathcal G(V(a))
[\mathcal V_E(V(a))-Y_k(V(a))]\,da.
\end{aligned}
\end{equation}
The integral is oriented when $t<0$, and the common trajectory preserves
\[
A_z(V(a))=A_z(U_0),\qquad
\mathscr S_R(V(a))=\mathscr S_R(U_0).
\]
\end{lemma}

\begin{proof}
We prove the claimed trajectory regularity and the chain rule \eqref{eq:legacy-53}, then integrate it to obtain \eqref{eq:legacy-54}.

For \(U_0\in\mathcal D\cap X^{11/2}\), set
\[
\gamma(a)=\Psi_k^{t-a}U_0,\qquad V(a)=\Phi_E^a\gamma(a),\qquad
\gamma\in C_a^1X^{11/2},\qquad\gamma'=-Y_k(\gamma).
\]
The original flow has a jointly \(C^1\) time/data map into \(X^{7/2}\). The two-derivative gap gives, for \(h\in X^{11/2}\),
\begin{equation}\label{eq:legacy-51}
\sup_{|t|\le T}\|\partial_tZ_h(t)\|_{X^{7/2}}
\le C_{M,T}\|h\|_{X^{11/2}},\qquad
\|D_U \Phi_E^{t+s}-D_U \Phi_E^t\|_{\mathcal L(X^{11/2},X^{7/2})}
\le C_{M,T}|s|.
\end{equation}
Indeed, \(\mathscr L_E:X^{11/2}\to X^{7/2}\) and the high variation bound imply \eqref{eq:legacy-51}. The \(C^2\) data map gives operator-norm continuity in the datum. These bounds establish a jointly \(C^1\) time/data map into \(X^{7/2}\).

Thus $V\in C_aX^{11/2}\cap C_a^1X^{7/2}$, and\[
V'(a)=\mathcal V_E(V(a))-D_U \Phi_E^a(\gamma(a))Y_k(\gamma(a)).
\]
The transport identity \eqref{eq:legacy-49} gives \eqref{eq:legacy-52}. The reference analytic map $\mathcal G$ acts on \(X\to Z_G\), and $X^{7/2}$ embeds in $X$. The strong chain rule therefore gives \eqref{eq:legacy-53}.
On these high data, the smooth kernel formula equals \(D_U \mathcal G\,\mathcal V_E\). The identity also holds in \(W_G\); only the weak difference is estimated uniformly in Section~\ref{sec:green}.

At the endpoints, \(V(0)=\Psi _k^tU_{0}\) and \(V(t)=\Phi _E^tU_{0}\). Integration of \eqref{eq:legacy-53} gives \eqref{eq:legacy-54}.
The integral is oriented for \(t<0\). The common high trajectories exist globally and preserve
\[
A_z(V(a))=A_z(U_0),\qquad \mathscr S_R(V(a))=\mathscr S_R(U_0).
\]
Equations \eqref{eq:legacy-52}, \eqref{eq:legacy-53}, and \eqref{eq:legacy-54}, together with these conservation identities, prove the lemma.
\end{proof}

\begin{remark}\label{rem:high-regularity-and-uniformity}
The high regularity in Lemma~\ref{lem:chain-rule} justifies differentiation of each fixed-parameter common trajectory. The uniform comparison in Lemma~\ref{lem:difference} uses only the endpoint bound and the common spectral tails; it does not require a bound for these high norms independent of $k$.
\end{remark}

\section{Green estimates and the real-line global flow}\label{sec:green}
We prove a scale-uniform weak bound and a strong approximation on fixed frequency bands. Conserved tails extend the comparison to actual common trajectories; the paired Green inverse recovers convergence in \(X\).

\subsection{Uniform Green bounds and fixed-band approximation}\label{subsec:uniform-green-bounds-and-fixed-band-approximation}

Two complementary estimates are needed: a weak bound uniform in the spectral parameter for rough backgrounds, and strong convergence for backgrounds with fixed Fourier support. Together they give estimates independent of the higher norms of the approximating trajectories.

Let $\mathcal G$ be the reference map \eqref{eq:reference-green-map}. Lemma~\ref{lem:determinant} makes this an analytic map \(X\to Z_G\), defined without pointwise rough kernels. Write \(\mathcal G(U)=(r_\eta,s_\eta)\in Z_G\) for the strong reference coordinates and \(\mathfrak V_E=D_U \mathcal G\,\mathcal V_E\) for the weak original field on smooth data. Their evolution is\begin{equation}\label{eq:legacy-55}
\mathfrak V_E(U)=(\sqrt3\partial_x s_\eta,(1-\partial_x^2-2u)\partial_x r_\eta/\sqrt3).
\end{equation}
For this calculation put $(\rho,\sigma,\tau)=(\rho_\eta,\sigma_\eta,\tau_\eta)$. The two kernel diagonals in $R_t=[\mathcal P_E,R]$, with $\mathcal P_E$ as in Lemma~\ref{lem:conservation}, verify this formula. The first component is \(\sqrt{3}\sigma '\). The second is initially \(\sqrt{3}\tau '+2u'\rho /\sqrt{3}\); use \(\rho '''+3\tau '=-4(\mathfrak a \rho )'\) to eliminate \(u'\rho \), leaving \((1-\partial_x^2-2u)\rho '/\sqrt{3}\). Derivatives of free constants vanish.

We estimate the transformed fields rather than the original endpoint vector field. The two-parameter formula gives a Lipschitz constant independent of the large spectral parameter in the weaker Green space.

\begin{lemma}\label{lem:UL}
On the ball of radius \(\varepsilon \), for \(k\ge 1\), the field $\mathfrak V_k(U):=D_U \mathcal G(U)Y_k(U)$ is well defined and\begin{equation}\label{eq:legacy-56}
\|\mathfrak V_k(U)-\mathfrak V_k(V)\|_{W_G}\le C_{\beta_*,\eta,\varepsilon}\|U-V\|_X.
\end{equation}
The original field in \eqref{eq:legacy-55} is also Lipschitz as a map \(X\to W_G\). In particular, the corresponding tail estimate is \(\|\mathfrak V_k(U)-\mathfrak V_k(P_{\le N}U)\|_{W_G}\le C\|P_{>N}U\|_X\), with a constant independent of $N,k$.
\end{lemma}

\begin{proof}
We prove \eqref{eq:legacy-56} for $k\ge2$ using the complete two-parameter formula, identify that formula with $D_U \mathcal G(U)Y_k(U)$ on rough backgrounds, and then cover $1\le k\le2$. The same product estimates give the original-field and tail bounds.

Estimate \eqref{eq:legacy-8} and \(C_k/|\zeta-\eta|\le Ck^2\) give \(p_k-p_k^0\in H^{3/2},q_k-q_k^0\in H^{1/2}\). Their norms and differences are controlled by $X$, uniformly in $k$.
The same bounds hold for $(\rho,\sigma)=(\rho_\eta,\sigma_\eta)$ after subtracting $(\rho_{0,\eta},\sigma_{0,\eta})$. At rough regularity these full diagonals mean the free constants plus $(r_\eta,s_\eta)$, rather than pointwise kernel values.
Apply Lemma~\ref{lem:products} to the products in \eqref{eq:legacy-35}.

Since $\mathfrak a=-\tfrac14+\tfrac12u$, we separate the constant part when estimating $\mathfrak a p$ and $\mathfrak a\rho$. The variable coefficient belongs to $H^{-1/2}$ and its derivative to $H^{-3/2}$. Use the last weak product for \(p''\rho',p'\rho'',q'\sigma,q\sigma'\) and the \(H^{3/2}\) multiplier bound for \(p t_\eta,t_k\rho\). Free constants act as identity multipliers and are not assigned \(L^2\) norms.
The first component lies in \(H^{-1/2}\) and the second in \(H^{-3/2}\).
Taking one difference in each term gives\begin{equation}\label{eq:legacy-57}
\|\mathcal F_k(U)-\mathcal F_k(V)\|_{W_G}
\le C_{\varepsilon,\beta_*,\eta}\|U-V\|_X,\quad k\ge2.
\end{equation}
Here \(\mathcal F_k\) temporarily denotes the rough Sobolev expression on the right of \eqref{eq:legacy-35}.

For fixed \(k\), choose real Schwartz \(U_j\to U\) inside a common spectral ball. The smooth identity and the preceding bounds give
\[
\mathcal F_k(U_j)=D_U \mathcal G(U_j)Y_k(U_j),\qquad
\mathcal F_k(U_j)\to\mathcal F_k(U)\text{ in }W_G,
\]
\[
D_U \mathcal G(U_j)Y_k(U_j)\to D_U \mathcal G(U)Y_k(U)\text{ in }Z_G.
\]
Since \(Z_G\hookrightarrow W_G\), the limits identify the original rough field. For \(1\le k\le2\), use
\[
D_U \mathfrak V_k[h]=D_U^2 \mathcal G[h,Y_k]+D_U \mathcal G\,D_U Y_k[h]
\]
and the compact-parameter field bounds. This covers equal spectral points. Applying \eqref{eq:legacy-56} to \(U,P_{\le N}U\) gives the tail estimate; \eqref{eq:legacy-55} obeys the same product bounds.
Thus \eqref{eq:legacy-56}, the original-field bound, and the tail estimate hold on the entire stated parameter range. This completes the proof.
\end{proof}

\begin{remark}\label{rem:rough-field-identification}
The density argument identifies the Sobolev formula \eqref{eq:legacy-35} with the original map $D_U \mathcal G(U)Y_k(U)$. Both maps are $C^1$ on the spectral ball, so their identity may also be differentiated there. This observation is used for the rough directions in the full-field comparison.
\end{remark}

For backgrounds supported in a fixed frequency band, the full auxiliary field approaches the original field in the strong coordinate norm. The remainder estimate below includes all differentiated trace orders and arbitrary rough test directions.

\begin{lemma}\label{lem:BA}
Fix \(N\ge 1\). On the ball of radius \(\varepsilon \),\begin{equation}\label{eq:legacy-58}
\sup_U\|\mathfrak V_k(P_{\le N}U)-\mathfrak V_E(P_{\le N}U)\|_{Z_G}\le C_{N,\beta_*,\varepsilon}/k.
\end{equation}
No spatial decay is required. Constants may depend on $N$, which is held fixed as $k$ tends to infinity.
\end{lemma}

\begin{proof}
We prove \eqref{eq:legacy-58} by estimating the quadratic, cubic, and complete higher-order parts of $Y_k-\mathcal V_E$ on a fixed band, and then apply $D_U \mathcal G$.

Fix $N\ge1$, $0<r\le\varepsilon$, and $S=P_{\le N}U$ with $\|U\|_X\le r$, with arbitrary direction \(h=(h_u,h_v)\in X\). For each fixed $m$, $\|S\|_{H^m\times H^{m-1}}+\|S\|_{L^\infty\times L^\infty}\le C_{N,m}r$.
The band-limited Hilbert--Schmidt estimate and vertex factorization \(TM_u\partial_x T+\partial_x TM_uT,\ TM_vT\) give\begin{equation}\label{eq:legacy-59}
\|B(S)\|_{HS}\le C_N r k^{-3/2},\quad
\|B(S)\|_{op}\le C_N r k^{-2},\quad
\|B(h)\|_{HS}\le Kk^{-1}\|h\|_X.
\end{equation}
The first two bounds may first be proved for \(k\ge 8N\). Increase the constants on a fixed bounded $k$ interval using the spectral ball.
The operator bound uses \(\|T\|\le Ck^{-3/2},\|\partial_x T\|\le Ck^{-1/2}\) and requires no extra decay of $u'$.

\smallskip\noindent\emph{Cubic test support.} The cubic variation is
\[
dA_z^{[3]}(S)[h]=\operatorname{Tr}(B_z(S)^2B_z(h)).
\]
Let the cyclic Fourier variables be $p_0,p_1,p_2$. Their differences \(\xi_1=p_0-p_1,\xi_2=p_1-p_2,\xi_3=p_2-p_0\) sum to zero.
The two $S$ vertices force \(|\xi_1|,|\xi_2|\le N\), hence \(|\xi_3|\le2N\).
The trace is therefore exactly unchanged when $h$ is replaced by \(P_{\le2N}h\).
The absolute cyclic integral is controlled by the Hilbert--Schmidt and bounded-operator ideal estimates for the absolute kernels. Apply Fubini before projecting.
Hilbert--Schmidt continuity extends the identity to rough $h$.
This support statement follows from the closed cycle of frequency differences and is used only at cubic order.

For \(k\ge8N\), put \(w_k(p)=(k^3+|p|^3)^{-1}\). Under shifts of size \(O(N)\),
\[
|m(p)|\lesssim w_k(p),\qquad |m'(p)|\lesssim k^{-1}w_k(p).
\]
Freezing the three \(u\)-vertices gives a multiplier error bounded by
\[
C_Nk^{-1}\int |p|^3w_k(p)^3\,dp\lesssim_N k^{-6}.
\]
Vertex errors and terms containing \(v\) have degree at most two, and
\[
\int(1+|p|)^2w_k(p)^3\,dp\lesssim k^{-6}.
\]

\smallskip\noindent\emph{Freezing both the background and spectral point.} We must include the displacement from \(z_k\) to \(k^3\) in the frozen integral. Set \(m_k^{\rm fr}(p)=(-ip^3-k^3)^{-1}\). The exact difference is\begin{equation}\label{eq:legacy-60}
m_{z_k}(p)-m_k^{\rm fr}(p)=\tfrac14(ip-k)m_{z_k}(p)m_k^{\rm fr}(p).
\end{equation}
For three $u$ vertices, the multiplier-difference integral is at most\[
C\int |p|^3(|p|+k)w_k(p)^4dp=Ck^{-7},
\]
The additional \(k\) term is still smaller than the required \(k^{-6}\) bound.
For the negative branch, the frozen multiplier is \((-ip^3+k^3)^{-1}\) and the exact difference is \((ip+k)m_{-z_k}m_{-,k}^{\rm fr}/4\), with the same absolute bound.
For the remaining frequency integral, Young's inequality gives
\[
\left|\iint f(\xi_1)g(\xi_2)h(-\xi_1-\xi_2)\,d\xi_1d\xi_2\right|
\le\|f\|_1\|g\|_2\|h\|_2.
\]
Band limitation bounds the two source factors by \(C_Nr\) and the projected test by \(C_N\|h\|_X\).

The frozen integral is\[
(2\pi)^{-1}\int(ip)^3[(ip)^3-k^3]^{-3}dp=1/(27k^5).
\]
The clockwise lower-half-plane residue calculation gives this sign. The negative branch follows by $p\mapsto-p$ and has the same value.
The three $u$ vertices contribute $1/8$ and the trace logarithm contributes $1/3$. Thus\begin{equation}\label{eq:legacy-61}
A_{\pm z_k}^{[3]}(S)=\frac{\int S_u^3}{648k^5}+\mathcal E_{3,\pm},\qquad
|d\mathcal E_{3,\pm}(S)[h]|\le C_Nr^2 k^{-6}\|h\|_X.
\end{equation}
The normalized average has leading term \(-\vartheta_k\int S_u^3/3\), with \(\vartheta_k=C_k/(216k^5)=(1-1/(12k^2))^3=1+O(k^{-2})\).
This compares the energy terms only for band-limited sources; no energy is defined for general rough data.

\smallskip\noindent\emph{Terms of order four and higher.} For \(n\ge 4\), retain the order of \(\operatorname{Tr}(B(S)^{n-1}B(h))\) and estimate it as follows:
use one Hilbert--Schmidt source factor, one Hilbert--Schmidt direction factor, and \(n-2\) operator-norm source factors. Hence\begin{equation}\label{eq:legacy-62}
\begin{split}
\sum_{n\ge4}\|B(S)^{n-1}B(h)\|_{\mathfrak S_1}
&\le\frac{\|B(S)\|_{HS}\|B(h)\|_{HS}\|B(S)\|_{op}^2}
{1-\|B(S)\|_{op}}\\
&\le C_{N,r,\beta_*}r^3 k^{-13/2}\|h\|_X.
\end{split}
\end{equation}
The denominator is bounded below by \(1-\beta_*\). The factor \(n\) from differentiation cancels \(1/n\); geometric convergence in trace norm justifies the entire differentiated series. At order \(n\), the test may reach \((n-1)N\), so the direction in \eqref{eq:legacy-62} remains arbitrary. After normalization,
\[
\|\nabla H_k^{[\ge4]}(S)\|_{X^*}\lesssim_{N,r}k^{-3/2},\qquad
\|Y_k^{[\ge4]}(S)\|_{X^{1/2}}\lesssim_{N,r}k^{-3/2}.
\]

\smallskip\noindent\emph{The quadratic term and the full field.} The exactly normalized quadratic field is $a_k(-i\partial_x)\mathscr L_ES$, where\[
a_k(\xi)=d_k^2/\mathfrak q_k(\xi).
\]
On the fixed band, \(|1-a_k|\le C_Nk^{-2}\).
The quadratic field error is at most $C_Nrk^{-2}$ in $X^{1/2}$.
Combining \eqref{eq:legacy-61} and \eqref{eq:legacy-62},\begin{equation}\label{eq:legacy-63}
\|Y_k(S)-\mathcal V_E(S)\|_{X^{1/2}}
\le C_N(rk^{-2}+r^2k^{-1}+r^3k^{-3/2}).
\end{equation}
The cubic normalization error is \(O(k^{-2})\) on the fixed band. Compact \(k\)-interval bounds extend \eqref{eq:legacy-63} to \(C_{N,r}/k\) for all \(k\ge1\). Since
\[
\|D_U \mathcal G(S)\|_{\mathcal L(X,Z_G)}\lesssim1,\qquad X^{1/2}\hookrightarrow X,
\]
the full Green field has the same approximation rate. Smooth approximation in a slightly larger fixed band identifies
\[
D_U \mathcal G(S)\mathcal V_E(S)=\mathfrak V_E(S)
\]
without spatial weights.

It follows that\begin{equation}\label{eq:legacy-64}
\sup_{\|U\|_X\le r}
\|\mathfrak V_k(P_{\le N}U)-\mathfrak V_E(P_{\le N}U)\|_{Z_G}\le C_{N,r,\beta_*,\eta}/k .
\end{equation}
Taking $r=\varepsilon$ in \eqref{eq:legacy-64} proves \eqref{eq:legacy-58} for every $k\ge1$.
This completes the proof.
\end{proof}
\begin{remark}\label{rem:band-hamiltonian-expansion}
The fixed-band estimate uses no energy or spatial-weight bound on $U$. The averaged quadratic and cubic terms converge to the corresponding energy terms on fixed bands. The single-branch translation term
\[
-\frac1{72\sqrt3 k^4}\int uv\,dx
\]
cancels in the average. The Hamiltonian $H_k$ itself is defined on $X$ by the determinant.
\end{remark}

\subsection{Strong recovery and the difference integral}\label{subsec:strong-recovery-and-the-difference-integral}

To recover the two physical fields we need a local inverse for the paired Green map. Its linear symbol is elliptic with the two Sobolev weights, and a small nonlinear perturbation preserves the two-point estimate.

\begin{lemma}\label{lem:inverse}
The reference linearization $M_1=D_U \mathcal G(0):X\to Z_G$ has symbol
\[
M_1(\xi)=-\frac1{2d_1\mathfrak q_1(\xi)}
\begin{pmatrix}\chi_1&\sqrt3\\1+\xi^2&\chi_1/\sqrt3\end{pmatrix}.
\]
It is an isomorphism. After decreasing \(\varepsilon \),\begin{equation}\label{eq:legacy-65}
\|U-V\|_X\le2\|M_1^{-1}\|\|\mathcal G(U)-\mathcal G(V)\|_{Z_G}.
\end{equation}
The image under $\mathcal G$ of any subset of the ball with a uniform $X$ Fourier tail has a uniform \(Z_G\) Fourier tail.
\end{lemma}

\begin{proof}
We compute the inverse of the reference linearization and absorb the nonlinear derivative error to prove \eqref{eq:legacy-65}. We then transfer a common Fourier tail through $\mathcal G$.

The single-branch quadratic gradient from Lemma~\ref{lem:coercivity} gives the stated symbol of $M_1$.
The determinant identity and inverse are
\[
\chi_1(\xi)^2-3(1+\xi^2)=\mathfrak q_1(\xi),\qquad
\begin{cases}
u=-2d_1[(13/4-\partial_x^2)r-3s],\\
v=2\sqrt3d_1[(1-\partial_x^2)r-(13/4-\partial_x^2)s].
\end{cases}
\]
Sobolev weights give the two bounded maps. On a smaller ball,
\[
\|M_1^{-1}(D_U \mathcal G(U)-M_1)\|_{\mathcal L(X)}\le\tfrac12.
\]
Integration along the segment from \(U\) to \(V\) and absorption prove \eqref{eq:legacy-65}.

For a subset $Q$ of the ball with common tail $\theta(N)=\sup_{U\in Q}\|P_{>N}U\|_X\to0$,
\[
\|\tau_yU-U\|_X\le CN|y|\varepsilon+2\theta(N),\qquad
\|\tau_y\mathcal G(U)-\mathcal G(U)\|_{Z_G}
\lesssim CN|y|\varepsilon+2\theta(N).
\]
Averaging the squared translation difference over \(|y|\le L^{-1}\) gives the multiplier
\[
2[1-\operatorname{sinc}(\xi/L)]\gtrsim1\qquad(|\xi|>L).
\]
Here $\operatorname{sinc}(r)=\sin r/r$, with value one at $r=0$. Consequently $\sup\|P_{>L}\mathcal G(U)\|_{Z_G}\le\Theta(L)\to0$.
Together with \eqref{eq:legacy-65}, this proves both assertions.
\end{proof}

For precompact \(Q\subset\mathcal D\), its closure remains in the larger analytic ball, even if it meets the strict spectral boundary. The comparison below uses high data in \(Q\) and requires no uniform sublevel gap.

The weak field estimate and the fixed-band approximation can now be combined along the actual common trajectories. The frequency cutoffs are fixed before the spectral parameter is sent to infinity.

\begin{lemma}\label{lem:difference}
For every precompact \(Q\subset \mathcal D \) and finite $T$,\begin{equation}\label{eq:legacy-66}
E_k(T,Q)=\sup_{U_0\in Q\cap X^{11/2},\,|t|\le T}\|\Phi_E^tU_0-\Psi_k^tU_0\|_X\longrightarrow0.
\end{equation}
\end{lemma}

\begin{proof}
We prove \eqref{eq:legacy-66} by combining the common spectral tails with Lemmas~\ref{lem:UL} and~\ref{lem:BA} in the weak Green norm, and then recover the strong endpoint norm using \eqref{eq:legacy-65}.

For \(U_0\in Q\cap X^{11/2}\), let
\[
V(a)=\Phi_E^a\Psi_k^{t-a}U_0,\qquad a\text{ between }0\text{ and }t.
\]
Lemma~\ref{lem:chain-rule} gives \(V\in C_aX^{11/2}\cap C_a^1X^{7/2}\) and
\begin{equation}\label{eq:legacy-67}
\partial_a\mathcal G(V)=(\mathfrak V_E-\mathfrak V_k)(V)
\quad\text{in }Z_G.
\end{equation}
Conservation and Lemmas~\ref{lem:coercivity} and~\ref{lem:spectral-tail-continuity} imply, uniformly in \(k,t,a\),
\begin{equation}\label{eq:legacy-68}
\|P_{>N}V(a)\|_X\le\theta_Q(N)
:=\left[c_2^{-1}\left(\sup_Q|\mathscr S_N|+C_3\varepsilon^3/N\right)\right]^{1/2}
\longrightarrow0.
\end{equation}
Lemmas~\ref{lem:UL} and~\ref{lem:BA} therefore give
\[
\|(\mathfrak V_E-\mathfrak V_k)(V)\|_{W_G}
\le C\theta_Q(N)+C_N/k.
\]
Integrating \eqref{eq:legacy-67},
\[
\|\mathcal G(\Phi_E^tU_0)-\mathcal G(\Psi_k^tU_0)\|_{W_G}
\le T(C\theta_Q(N)+C_N/k).
\]
Lemma~\ref{lem:inverse} supplies common strong tails \(\Theta_Q(L)\to0\). Since the strong and weak coordinates differ by two derivatives, \eqref{eq:legacy-65} yields
\begin{equation}\label{eq:legacy-69}
E_k(T,Q)\le C\|M_1^{-1}\|
\left[\langle L\rangle^2T(C\theta_Q(N)+C_N/k)+2\Theta_Q(L)\right].
\end{equation}
Choose \(L\), then \(N\), then \(k\). Equivalently,
\[
\lim_{L\to\infty}\lim_{N\to\infty}\limsup_{k\to\infty}E_k(T,Q)=0.
\]
This proves \eqref{eq:legacy-66}, with no dependence on the high norms of the approximations.
This completes the proof.
\end{proof}

\subsection{Construction of the real-line global group}\label{subsec:construction-of-the-real-line-global-group}

We pass from the uniform smooth comparison to the endpoint global group. The spectral parameter is chosen before invoking fixed-parameter continuity; its growing stability constant therefore causes no loss. The local equation and uniqueness class are identified in Section~\ref{subsec:local-identification-and-the-canonical-source}.

\begin{lemma}\label{lem:global-limit}
For $U_0\in\mathcal D$ and any real Schwartz sequence $U_{0,n}\in\mathcal D$ with $U_{0,n}\to U_0$ in $X$, the limit
\[
\Phi^\cdot U_0=\lim_{n\to\infty}\Phi_E^\cdot U_{0,n}
\quad\text{in }C([-T,T];X)
\]
exists for every $T>0$ and is independent of the sequence. It defines the continuous global group and spectral conservation asserted in Theorem~\ref{thm:line}.
\end{lemma}

\begin{proof}
We prove that the high-regularity approximations are Cauchy in $C([-T,T];X)$, establish independence and spectral conservation, and then verify data continuity and the group law.

Choose real Schwartz \(U_n\in\mathcal D\) with \(U_n\to U_0\) in \(X\). For \(Q=\{U_0,U_n:n\ge1\}\), \eqref{eq:legacy-66} and fixed-\(k\) continuity imply
\begin{equation}\label{eq:legacy-70}
\|\Phi_E^\cdot U_n-\Phi_E^\cdot U_m\|_{C_tX}
\le2E_k(T,Q)+K_{k,T}\|U_n-U_m\|_X,\qquad K_{k,T}=e^{Ck^3T}.
\end{equation}

First choose \(k\), then \(n,m\). This defines
\[
\Phi^\cdot U_0=\lim_{n\to\infty}\Phi_E^\cdot U_n\quad\text{in }C([-T,T];X).
\]
Interleaving approximating sequences proves independence. Spectral continuity gives
\[
A_z(\Phi^tU_0)=A_z(U_0),\qquad\mathscr S_R(\Phi^tU_0)=\mathscr S_R(U_0),\qquad
\|\Phi^tU_0\|_X<\varepsilon/2.
\]

The same comparison applies to any sequence $V_n\in\mathcal D\cap X^{11/2}$ converging to $U_0$ in $X$. Indeed, the union of this sequence and a Schwartz approximation is precompact, and \eqref{eq:legacy-66} compares each of its high-regularity trajectories with the same fixed-$k$ auxiliary flow. Choosing $k$ before $n$ in \eqref{eq:legacy-70} therefore gives
\[
\Phi_E^\cdot V_n\longrightarrow\Phi^\cdot U_0
\quad\text{in }C([-T,T];X).
\]

For continuous dependence, let \(U_j\to U\in\mathcal D\). Choose \(\rho_j>0\) with \(B_X(U_j,\rho_j)\subset\mathcal D\), and real Schwartz approximations
\begin{equation}\label{eq:legacy-71}
S_{j,n}\in\mathcal D,\qquad
\|S_{j,n}-U_j\|_X<\min\{(j+n)^{-1},\rho_j/2\}.
\end{equation}
Include a sequence for \(U\). The union is totally bounded: large \(j\) places it near \(U\), and each of the finitely many remaining rows is convergent. Its only accumulation points are $U$ and the $U_j$, all in $\mathcal D$. Hence the compact closure $K$ lies in $\mathcal D$. Passing to \(n\to\infty\) in \eqref{eq:legacy-70},
\[
\|\Phi^\cdot U_j-\Phi^\cdot U\|_{C_tX}
\le2E_k(T,K)+K_{k,T}\|U_j-U\|_X\longrightarrow0
\]
by choosing \(k\) before \(j\). Time continuity follows from the strong limit. At intermediate times,
\[
\Phi^t\Phi^sU=\lim_n\Phi_E^t\Phi_E^sU_n=\Phi^{t+s}U,\qquad
\Phi^0=\mathrm{id},\qquad(\Phi^t)^{-1}=\Phi^{-t}.
\]
This proves the approximation, conservation, continuity, and group assertions.
\end{proof}

\begin{remark}\label{rem:limiting-auxiliary-commutation}
For fixed $k$ and $b$, the actual commutation of the high-regularity flows, their endpoint convergence, and the continuity of $\Psi_k^b$ also give
\[
\Phi^t\Psi_k^bU=\Psi_k^b\Phi^tU,\qquad U\in\mathcal D.
\]
This additional identity is not needed to define the group in Lemma~\ref{lem:global-limit}.
\end{remark}

\subsection{Local identification and the canonical source}\label{subsec:local-identification-and-the-canonical-source}

Strong convergence of the trajectories does not by itself identify their nonlinear sources. We use the local bilinear source estimate and low-frequency recovery to identify the equation and its prescribed uniqueness class.

\begin{lemma}\label{lem:source}
The limit in Lemma~\ref{lem:global-limit} has a canonical differentiated source
\[
\mathcal N(U)=\partial_x (u^2)_{\rm can}\in\mathcal D'_tH^{-7/2}_x,
\qquad u_t=\partial_x v,\quad v_t=\partial_x (u-\partial_x^2u)-\mathcal N(U).
\]
The source is independent of the high-regularity approximation, and $\partial_x \mathcal N$ equals the source in the local second-order equation. The solution is unique in the fixed class $\mathcal K_{\lambda_*}$.
\end{lemma}

\begin{proof}
We identify the local differentiated source, prove its independence of extensions and approximations, and recover the two-field equation. We then establish uniqueness in the fixed class $\mathcal K_{\lambda_*}$.

\emph{Source estimates.} Apply the linear and bilinear estimates recalled in Section~\ref{subsec:local}, with \(\lambda=\lambda_*\). In rescaled variables \((t_{\mathrm{loc}},y)\), write
\[
\Lambda=\langle\omega+q^2\rangle,\qquad
\Omega_\lambda(q)=\frac{\lambda^2q^2}{1+\lambda^2q^2},\qquad
\mathcal F_\lambda=\Lambda W,\quad
\|F\|_{\mathcal F_\lambda}=\|\Lambda^{-1}F\|_W.
\]
The embedding \(W\hookrightarrow X_{\mathrm{Sch}}^{-1/2,0}\) gives
\begin{equation}\label{eq:legacy-72}
|\langle F,\chi\rangle|
\le C_p\|F\|_{\mathcal F_\lambda}
\|\langle q\rangle^{1/2}\langle\omega+q^2\rangle\widehat\chi\|_2.
\end{equation}
Thus \(\mathcal F_\lambda\hookrightarrow\mathcal S'\). Define by smooth density
\[
Q_\lambda(b,c)=\tfrac14\Omega_\lambda(b+\bar b)(c+\bar c).
\]
The bilinear and linear local estimates yield
\begin{equation}\label{eq:legacy-73}
\|Q_\lambda(b,c)\|_{\mathcal F_\lambda}
\le C_p\lambda^{1/2}\|b\|_W\|c\|_W,
\end{equation}
\begin{equation}\label{eq:legacy-74}
\begin{gathered}
F_\lambda(b)=\tfrac12\lambda^{-2}(b-\bar b)-Q_\lambda(b,b),\\
\|F_\lambda(b)-F_\lambda(c)\|_{\mathcal F_\lambda}
\le C_p\bigl[\lambda^{-2}+\lambda^{1/2}(\|b\|_W+\|c\|_W)\bigr]\|b-c\|_W.
\end{gathered}
\end{equation}

\emph{Locality.} For a compactly supported unit-integral approximate identity \(\rho_\delta\), dominated convergence in each Fourier component, followed by the dyadic \(\ell^1\) tail estimate, gives
\[
\rho_\delta*b\longrightarrow b\quad\text{in }W.
\]
The mollified functions have all spacetime derivatives in \(L^2\). Schwartz approximation in these norms and \(W\) identifies \(Q_\lambda\) with their ordinary smooth product. If \(b=b'\), \(c=c'\) on an open time interval \(I\), then, for \(I'\Subset I\) and sufficiently small \(\delta\),
\[
\rho_\delta*b=\rho_\delta*b',\qquad\rho_\delta*c=\rho_\delta*c'
\quad\text{on }I'\times\mathbb R.
\]
Since \(\Omega_\lambda\) is spatial, \eqref{eq:legacy-72}--\eqref{eq:legacy-74} imply
\[
Q_\lambda(b,c)|_I=Q_\lambda(b',c')|_I,\qquad
F_\lambda(b)|_I=F_\lambda(b')|_I.
\]
Restriction source norms are therefore well defined. Nearly optimal extensions give the corresponding restriction estimates. On \(|t_{\mathrm{loc}}|<1\),
\[
(i\partial_{t_{\mathrm{loc}}}+\partial_y^2)\mathcal IF=F,
\]
first for smooth sources and then by the extension estimate and \eqref{eq:legacy-72}.

\emph{Local identification.} Fix \(t_0\). For smooth original trajectories \(U_n\to U\) in \(C_tX\), Lemma~\ref{lem:persistence} supplies extensions satisfying
\[
b_n\to b\text{ in }W,\qquad F_\lambda(b_n)\to F_\lambda(b)\text{ in }\mathcal F_\lambda.
\]
High mild uniqueness identifies their restrictions with
\[
b_n(t_{\mathrm{loc}},y)=\lambda^{-2}w_n(t_0+\lambda^{-2}t_{\mathrm{loc}},y/\lambda).
\]
On \(I=t_0+\lambda^{-2}(-1,1)\), put
\begin{equation}\label{eq:legacy-75}
u(t,x)=\lambda^2\Re b(\lambda^2(t-t_0),\lambda x),\qquad
P=\lambda^4Q_\lambda(b,b)(\lambda^2(t-t_0),\lambda x)
=\Omega_1(u^2)_{\mathrm{loc}}.
\end{equation}
For spacetime tests,
\begin{equation}\label{eq:legacy-76}
\langle P,\chi\rangle
=\lambda\langle Q_\lambda(b,b),
\chi(t_0+t_{\mathrm{loc}}/\lambda^2,y/\lambda)\rangle.
\end{equation}
The fixed scaling and \eqref{eq:legacy-72}--\eqref{eq:legacy-74} give
\begin{equation}\label{eq:legacy-77}
\mathcal N_2(U):=-(1-\partial_x^2)P
=\lim_n \partial_x^2(u_n^2)\quad\text{in }\mathcal D'(I\times\mathbb R).
\end{equation}
This is also the consequence of
\[
\partial_y^2[(\Re b_n)^2]
=-\lambda^{-2}(1-\lambda^2\partial_y^2)Q_\lambda(b_n,b_n).
\]

\emph{The first-order source.} Strong convergence gives \(u_t=\partial_x v\). For a time test \(\phi\), define
\begin{equation}\label{eq:legacy-78}
\mathcal N_U[\phi]=\int_I\phi \partial_x (u-\partial_x^2u)\,dt+\int_I\phi'v\,dt
\in H^{-7/2}.
\end{equation}
Then
\[
\|\mathcal N_U[\phi]\|_{H^{-7/2}}
\lesssim\|U\|_{C_tX}(\|\phi\|_1+\|\phi'\|_1),\qquad
\mathcal N_{U_n}=\partial_x (u_n^2)\to\mathcal N_U.
\]
Comparing this limit with \eqref{eq:legacy-77} yields
\begin{equation}\label{eq:legacy-79}
\partial_x \mathcal N_U=\mathcal N_2(U)\quad\text{in }\mathcal D'(I;H^{-9/2}).
\end{equation}
For every \(s\in\mathbb R\),
\[
f\in H^s,\quad \partial_x f=0
\quad\Longrightarrow\quad\operatorname{supp}\widehat f\subset\{0\}
\quad\Longrightarrow\quad f=0.
\]
The same holds after pairing any Sobolev-valued time distribution with a time test. Therefore \(\mathcal N_U\) is the unique \(H^{-7/2}\)-valued first-order source satisfying \eqref{eq:legacy-79}. Write
\[
\partial_x (u^2)_{\mathrm{can}}=\mathcal N_U,\qquad
v_t=\partial_x (u-\partial_x^2u)-\mathcal N_U.
\]
Moreover,
\[
u_t=(1-\partial_x^2)\Im w\in C_tH^{-5/2},\qquad
\partial_x v(t_0)=u_t(t_0),\qquad v(t_0)\in H^{-3/2}.
\]
The zero-kernel argument fixes \(v\) uniquely without a low-frequency \(\partial_x^{-1}\).

\emph{Independence and uniqueness.} Interleaving the smooth sequences gives the same \(U\); continuity of the fixed local map gives the same source. Locality identifies different extensions, and \eqref{eq:legacy-78} identifies overlapping windows.

For two integral solutions \(b,c\) in \(\mathcal K_{\lambda_*}\) with equal initial data, the short-interval zero-data estimate in \cite[Section 2]{r1} gives
\[
\|b-b_{\mathrm{free}}\|_{W(I)}+
\|c-b_{\mathrm{free}}\|_{W(I)}\longrightarrow0
\quad(|I|\to0).
\]
Choose the initial small ball and \(\lambda_*\) as in Lemma~\ref{lem:persistence}, then shorten \(I\). The restriction extension and bilinear estimates yield
\[
\|b-c\|_{W(I)}
\le C_p\bigl[\lambda_*^{-2}+\lambda_*^{1/2}(\|b\|_{W(I)}+\|c\|_{W(I)})\bigr]
\|b-c\|_{W(I)}
\le\tfrac12\|b-c\|_{W(I)}.
\]
Thus \(b=c\). Iteration gives equality for both time directions; \(u_t=\partial_x v\) and the zero-kernel argument identify \(v\). All comparisons use the same fixed class.
The canonical source, its compatibility with the local equation, and fixed-class uniqueness have now been proved. This completes the proof.
\end{proof}

\begin{proof}[Proof of Theorem~\ref{thm:line}]
We verify the three assertions using the preceding construction. Lemma~\ref{lem:global-limit} gives the global group on $\mathcal D$, continuity into $C([-T,T];X)$, spectral conservation, and convergence for every approximating sequence in $\mathcal D\cap X^{11/2}$. In particular,
\[
\sup_{t\in\mathbb R}\|\Phi^tU_0\|_X^2
\le\gamma^{-1}\mathscr S(U_0)<\varepsilon^2/4.
\]
For $U_j\to U$ in $\mathcal D$ and $t_j\to t$, choose $T>|t|$ with $|t_j|\le T$ for large $j$. Then
\[
\|\Phi^{t_j}U_j-\Phi^tU\|_X
\le\|\Phi^\cdot U_j-\Phi^\cdot U\|_{C([-T,T];X)}
+\|\Phi^{t_j}U-\Phi^tU\|_X\longrightarrow0.
\]
This proves joint continuity and all assertions in parts (1) and (2). Lemma~\ref{lem:source} gives the distributional equation with the canonical differentiated source and uniqueness in $\mathcal K_{\lambda_*}$, which is part (3). This completes the proof of Theorem~\ref{thm:line}.
\end{proof}

\section{The periodic problem}\label{sec:periodic}
The periodic result follows from the discrete spectral functional and a fixed endpoint local theory. We prove the spectral bounds on the full periodic operator space, then obtain a uniform local lifespan, persistence, and a canonical differentiated source. Iteration gives the global group.

\subsection{The discrete spectral functional}\label{subsec:the-discrete-spectral-functional}

We repeat the spectral coercivity calculation with the actual discrete Fourier sums. Low modes and the constant spectral mode remain present; the positivity argument does not rely on a continuum approximation.
Throughout this section, unadorned norms of two-field data mean $\|\cdot\|_{X_{\mathbb T}}$, and
\[
\mathbb T=\mathbb R/(2\pi\mathbb Z),\qquad
d\mu=\frac{dx}{2\pi},\qquad f_n=\int_{\mathbb T}f(x)e^{-inx}\,d\mu,
\]
\[
\|U\|_{X^s_{\mathbb T}}^2
=\sum_n\langle n\rangle^{2s}|u_n|^2
+\sum_n\langle n\rangle^{2s-2}|v_n|^2,
\qquad X_{\mathbb T}=X^{-1/2}_{\mathbb T}.
\]
Let $X_{\mathbb T,0}$ be its real subspace with $\int_{\mathbb T}u\,d\mu=\int_{\mathbb T}v\,d\mu=0$. The spectral operators act on the full $L^2(\mathbb T,d\mu)$, including the constant mode. Put
\[
a_n=n^3+n/4,\quad z_k=k^3-k/4,\quad
m_\pm(n)=(-ia_n\mp z_k)^{-1},\quad t_\pm(n)^2=m_\pm(n),\quad k\ge1.
\]
The square-root branches are those fixed in the introduction. In particular, $m_\pm(0)=\mp z_k^{-1}$. Define
\begin{equation}
(B_\pm(U))_{nm}
=t_\pm(n)t_\pm(m)
\left[\frac{i(n+m)}4u_{n-m}+\frac{v_{n-m}}{4\sqrt3}\right].
\label{eq:torusB}
\end{equation}

The periodic spectrum requires a discrete calculation at every mode. The next lemma keeps the constant spectral mode, cancels the mixed quadratic term by the two-branch average, and controls the entire trace tail.

\begin{lemma}\label{lem:periodic-spectral}
There are $r,c_2,C_2,C_3>0$ such that, for real $U\in X_{\mathbb T,0}$ with $\|U\|<r$,
\begin{equation}
\begin{gathered}
A^\mathbb T_\pm(k;U)
=\sum_{j\ge2}\frac{(-1)^{j+1}}j\operatorname{Tr}B_\pm(U)^j,
\qquad s^\mathbb T_k=-\tfrac12(A^\mathbb T_++A^\mathbb T_-),\\
\mathscr S^\mathbb T_R(U)=\int_R^\infty k s^\mathbb T_k(U)\,dk,
\qquad \mathscr S^\mathbb T=\mathscr S^\mathbb T_1,\\
\frac{c_2}{2}\|U\|^2\le\mathscr S^\mathbb T(U)
\le(C_2+C_3r)\|U\|^2.
\end{gathered}
\label{eq:toruscoercive}
\end{equation}
These are continuous real functionals. They are conserved along smooth periodic solutions for as long as those solutions remain in the spectral ball.
\end{lemma}

\begin{proof}
We prove convergence and reality of the trace series, compute the discrete quadratic term and control the full remainder to obtain \eqref{eq:toruscoercive}. We then justify conservation on existing smooth trajectories.

Set $b_n=(k+|n|)^{-3}$. Since $|m_\pm(n)|\lesssim b_n$,
\begin{equation}
\sum_n b_nb_{n+l}\lesssim k^{-2}(k+|l|)^{-3},\qquad
\sum_n(2n+l)^2b_nb_{n+l}\lesssim k^{-2}(k+|l|)^{-1}.
\label{eq:torusconvolution}
\end{equation}
For $|l|\le2k$, split at $|n|=4k$. The two sums on the inner region are $O(k^{-5})$ and $O(k^{-3})$; the outer tails satisfy the same bounds. For $|l|>2k$, use the regions $|n|\le|l|/2$, $|n+l|\le|l|/2$, and their complement. The first two use $\sum b_n\lesssim k^{-2}$; the remaining tails are $O(|l|^{-5})$ and $O(|l|^{-3})$. Thus
\begin{equation}
\|B_\pm(U)\|_{\mathfrak S_2}\le C_Hk^{-1}\|U\|_{X_{\mathbb T}}.
\label{eq:torusHS}
\end{equation}
Choose $C_Hr\le\beta<1$. The trace series converges absolutely and is the branch of $\log\det_2(I+B_\pm)$ vanishing at zero. Finite-rank approximation in $\mathfrak S_2$ proves the determinant identity; the inverse is the Neumann inverse. Fourier reflection shows that each trace is real. No self-adjointness is assumed.

For $l\ne0$, write
\[
J_j(k,l)=\sum_n(2n+l)^j m_+(n)m_+(n+l),\qquad
P_\pm=3n^2\pm3nl+l^2+\tfrac14.
\]
The sums are absolute. The identities
\[
\begin{gathered}
m_+(n)-m_+(n+l)=-ilP_+m_+(n)m_+(n+l),\\
(2n+l)^2=\tfrac43P_+-\tfrac13(l^2+1)
\end{gathered}
\]
and an integer shift give
\begin{equation}
J_2=-\tfrac13\langle l\rangle^2J_0,\qquad
J_0=6\sum_n\frac{n a_n}{(z_k^2+a_n^2)P_+P_-}>0.
\label{eq:torusJ}
\end{equation}
Indeed, $J_0=-6i\sum n m_+(n)/(P_+P_-)$, and the odd part vanishes. Since $P_\pm\sim n^2+l^2+1$,
\begin{equation}
J_0\sim\sum_{n\ge1}\frac{n^4}{(k^6+n^6)(n^2+l^2)^2}
\sim\frac1{k(k^2+l^2)^2}.
\label{eq:torusJbounds}
\end{equation}
For the upper bound, use $|l|\le k$ and $|l|>k$ separately, with
$\sum(k^6+n^6)^{-1}\lesssim k^{-5}$ and
$\sum n^4/(k^6+n^6)\lesssim k^{-1}$.
For the lower bound, retain $\lceil k\rceil\le n\le2\lceil k\rceil$.
All constants are uniform for $k\ge1$ and $l\ne0$.

Reflection gives $J_j^- =(-1)^jJ_j$. Consequently the mixed quadratic term cancels in the spectral average and
\begin{equation}
(s^\mathbb T_k)^{[2]}(U)
=\sum_{l\ne0}\frac{\langle l\rangle^2J_0(k,l)}{96}
\left(|u_l|^2+\frac{|v_l|^2}{\langle l\rangle^2}\right).
\label{eq:torusquadratic}
\end{equation}
The internal spectral sum in $J_0$ still includes $n=0$. Only the data mode $l=0$ is absent. Integration of \eqref{eq:torusJbounds} yields
\[
\int_1^\infty k\frac{\langle l\rangle^2J_0(k,l)}{96}\,dk
\sim\langle l\rangle^{-1}.
\]
For the full higher-order tail,
\begin{equation}
\begin{aligned}
|s^\mathbb T_k-(s^\mathbb T_k)^{[2]}|
&\le\sum_{j\ge3}\frac{(C_Hk^{-1}\|U\|)^j}{j}
\le C_3 k^{-3}\|U\|^3,\\
|d(s^\mathbb T_k-(s^\mathbb T_k)^{[2]})[H]|
&\le C'_3k^{-3}\|U\|^2\|H\|.
\end{aligned}
\label{eq:torusremainder}
\end{equation}
Tonelli applies to the positive quadratic term. The higher-order integral is absolute because $\int_R^\infty k^{-2}\,dk=R^{-1}$. Choosing $C_3r\le c_2/2$ gives \eqref{eq:toruscoercive} and continuity, with no endpoint logarithmic loss.

It remains to justify conservation. For a smooth solution, the periodic Lax equation is
\[
L_t=[P,L],\qquad P=P_0+M,\quad
P_0=\sqrt3\partial_x^2,\quad M=(u-\tfrac12)/\sqrt3.
\]
Fix either spectral sign and let $T$ be the multiplier with symbol $t_\pm(n)$,
$H_+=TMT^{-1}$, $H_-=T^{-1}MT$, and $C=T[M,L_0]T$.
Smooth Fourier decay and $|t(n)/t(m)|\lesssim_k\langle n-m\rangle^{3/2}$ give
\[
H_\pm\in\mathcal B(L^2),\qquad C=H_+-H_-\in\mathfrak S_2,
\qquad [P_0,B],B_t\in\mathfrak S_2.
\]
For example the matrix of $C$ is bounded by
$c(k)\langle l\rangle^{9/2}\langle n\rangle^{-1}|M_l|$, $l=n-m$, with $c(k)$ depending on the fixed spectral point.
Direct multiplication gives
\begin{equation}
B_t=[P_0,B]+[H_+,B]+(I+B)C.
\label{eq:torusBt}
\end{equation}
The unbounded $P_0$ is handled by its unitary group: differentiation of
$\operatorname{Tr}(e^{isP_0}Be^{-isP_0})^j$ in $\mathfrak S_2$ gives
$\operatorname{Tr}B^{j-1}[P_0,B]=0$. Cyclicity for $H_+$ is legitimate. Hence
\[
\partial_t\sum_{j=2}^J\frac{(-1)^{j+1}}j\operatorname{Tr}B^j
=-\operatorname{Tr}BC+(-1)^{J+1}\operatorname{Tr}B^JC.
\]
The last term tends to zero uniformly on compact smooth time intervals. In $\operatorname{Tr}BC$, the identity
\[
m(n)m(m)[-ia_m+ia_n]=m(n)-m(m)
\]
reduces the $v$ contribution to a shifted sum of $m$ and the $u$ contribution to
\[
\sum_m(2m+l)[m(m+l)-m(m)]=-2l\sum_m m(m).
\]
Here $m$ and $nm(n)$ are summable. The remaining factor is
$\sum_l il u_lu_{-l}=0$. Thus $A^\mathbb T_\pm$ is conserved. Integrating the pointwise spectral equality proves conservation of $\mathscr S^\mathbb T_R$ without differentiating under the spectral integral.
This proves all the asserted spectral bounds and conservation laws.
\end{proof}

\subsection{A fixed periodic local theory}\label{subsec:a-fixed-periodic-local-theory}

We use the variables $w$ and $b=b_{\lambda,0}$ from Remark~\ref{rmk:local-scales} and the periodic space $W=\mathcal W_{p,\lambda}$ from Definition~\ref{def:Wp}. The zero-mean condition gives the additional recovery formula
\[
v=\partial_x^{-1}(1-\partial_x^2)\Im w,
\]
where $\partial_x^{-1}$ is zero on the constant mode. Since $2^{-1/2}\le |n|/\langle n\rangle\le1$ for $n\ne0$,
\begin{equation}\label{eq:toruscoordinate}
2^{-1/2}\|U\|_{X^s_{\mathbb T}}
\le\|w\|_{H^s(\mathbb T)}\le\|U\|_{X^s_{\mathbb T}}.
\end{equation}
For the fixed scale $\lambda\ge1$ selected below, $b$ satisfies
\begin{equation}\label{eq:torusscaling}
ib_s+\partial_y^2b=\tfrac12\lambda^{-2}(b-\bar b)
-\Omega_\lambda(\Re b)^2.
\end{equation}
The multiplier $\Omega_\lambda$ and the circle $\mathbb T_\lambda$ are those fixed in Remark~\ref{rmk:local-scales}.
This changes both the period and a coefficient. It is not a symmetry of the fixed-period equation.

In $W$ we retain Kishimoto's spatial Fourier convention on $\mathbb T_\lambda$ and the frequency measure $\lambda^{-1}\sum_{q\in\lambda^{-1}\mathbb Z}$ from Definition~\ref{def:Wp}. The scaling of the initial datum is therefore
\[
\|b(0)\|_{H^{-1/2}(\mathbb T_\lambda)}^2
=2\pi\lambda^{-3}\sum_n\langle n/\lambda\rangle^{-1}|w_n(0)|^2,
\]
so the forward and inverse scaling norms are at most
\begin{equation}
C_f=\sqrt{2\pi}\lambda^{-1},\qquad
C_i=\lambda^{3/2}/\sqrt{2\pi}.
\label{eq:torusscalingbounds}
\end{equation}
Set $\Lambda=\langle\sigma+q^2\rangle$ and $\mathcal F_\lambda=\Lambda W$. The special extension $\mathcal I$ is the three-term extension of Section~\ref{sec:preliminaries}, with periodic spatial Fourier transform. Its last term is not multiplied by a time cutoff. Kishimoto's periodic estimates \cite{r1} yield
\begin{equation}
\begin{gathered}
\|\psi(s)e^{is\partial_y^2}f\|_W\le C_p\|f\|_{H^{-1/2}},\qquad
\|\mathcal IF\|_W\le C_p\|F\|_{\mathcal F_\lambda},\\
\|\tfrac12\lambda^{-2}\mathcal I(b-\bar b)\|_W
\le C_p\lambda^{-2}\|b\|_W,\\
\|Q_\lambda(b,c)\|_{\mathcal F_\lambda}
\le C_p\lambda^{1/2}\|b\|_W\|c\|_W,
\quad Q_\lambda(b,c)=\Omega_\lambda(\Re b\Re c),\\
W\hookrightarrow C_sH^{-1/2},\qquad
X_{\rm Sch}^{-1/4,3/4}\hookrightarrow W.
\end{gathered}
\label{eq:torus-local-estimates}
\end{equation}
The mixed source has the factor $\lambda^{1/2}$; the two same-sign products do not. These are source estimates, not only estimates for their Duhamel images. The modulation decomposition and Duhamel extension are those fixed in Definition~\ref{def:Wp} and Section~\ref{subsec:local}, with the periodic spatial Fourier transform. The short-window property needed for uniqueness is proved in Lemma~\ref{lem:zero-trace}.

For $f=b(0)$, define the linear maps $\mathcal A,\mathcal L$ and the symmetric real-bilinear map $\mathcal B$ by
\begin{equation}\label{eq:periodic-fixed-point-operators}
\begin{aligned}
(\mathcal Af)(s)&=\psi(s)e^{is\partial_y^2}f,\\
\mathcal Lb&=\tfrac12\lambda^{-2}\mathcal I(b-\bar b),\\
\mathcal B(b,c)&=-\mathcal I Q_\lambda(b,c)
=-\mathcal I\bigl[\Omega_\lambda(\Re b\,\Re c)\bigr].
\end{aligned}
\end{equation}
Here $\mathcal I$ is the periodic special Duhamel extension fixed above; on $|s|<1$ it equals $-i\int_0^s e^{i(s-a)\partial_y^2}(\cdot)(a)\,da$. Thus its factor $-i$ is already included in $\mathcal L$ and $\mathcal B$.

The bounded periodic recovery transfers local stability from the complex variable to both fields. The scale is selected once, so the following lifespan remains available at every small-data restart.

\begin{lemma}\label{lem:periodic-local}
There are fixed $\lambda,r_{\rm loc},\tau>0$ and $L_{\rm loc}\ge1$ with the following properties. For $U_0\in X_{\mathbb T,0}$ with $\|U_0\|\le r_{\rm loc}$, the rescaled initial datum $f=b(0)$ has a unique extension in the contraction ball satisfying
\[
b=\mathcal Af+\mathcal Lb+\mathcal B(b,b),
\]
with the operators \eqref{eq:periodic-fixed-point-operators}. The recovered solution is defined on $[-\tau,\tau]$, where $\tau=\lambda^{-2}$, and
\begin{equation}
\|U-V\|_{C([-\tau,\tau];X_{\mathbb T})}
\le L_{\rm loc}\|U(0)-V(0)\|_{X_{\mathbb T}}.
\label{eq:toruslocaldifference}
\end{equation}
If $U_0\in X^{j-1/2}_{\mathbb T}$, $j\ge0$ an integer, the same interval preserves that regularity. Neither $\tau$ nor the endpoint difference constant depends on a higher initial norm.
\end{lemma}

\begin{proof}
We choose one contraction scale and radius, derive \eqref{eq:toruslocaldifference}, and prove persistence on that same lifespan by differentiating the translation orbit.

By \eqref{eq:periodic-fixed-point-operators} and \eqref{eq:torus-local-estimates},
\[
\|\mathcal A\|\le C_p,\qquad
\|\mathcal L\|\le C_p\lambda^{-2},\qquad
\|\mathcal B\|\le C_p\lambda^{1/2}.
\]
Choose $\lambda$, then $r_{\rm loc}$, so that
\begin{equation}
C_p\lambda^{-2}\le\tfrac18,\qquad
R_b=2C_p C_f r_{\rm loc},\qquad
2C_p\lambda^{1/2}R_b\le\tfrac18.
\label{eq:toruscontraction}
\end{equation}
The map preserves $B_W(0,R_b)$ and contracts with constant $q\le1/4$. Subtraction gives
\[
(1-q)\|b-c\|_W\le C_p\|f-g\|_{H^{-1/2}}.
\]
Combining \eqref{eq:toruscoordinate}--\eqref{eq:torusscalingbounds} and the trace embedding gives \eqref{eq:toruslocaldifference}, with
\[
L_{\rm loc}=\max\{1,\sqrt2 C_iC_{\rm emb}C_p(1-q)^{-1}C_f\}.
\]
The fixed-point map is real analytic into $W$, since
\[
\|(I-\mathcal L-2\mathcal B(b,\cdot))^{-1}\|\le(1-q)^{-1}.
\]
Spatial translations are strongly continuous in $W$: apply dominated convergence to the weighted $L^2$ and $Y$ terms, and first discard the summable tail in the dyadic $\ell^1$ term. The solution map commutes with translations. Differentiating the initial translation orbit in $H^{-1/2}$ therefore identifies $\partial_y^jb\in W$ for $f\in H^{j-1/2}$. Equivalently,
\begin{equation}
(I-\mathcal L-2\mathcal B(b,\cdot))\partial_y^jb
=\mathcal A\partial_y^jf+\sum_{a=1}^{j-1}\binom ja\mathcal B(\partial_y^ab,\partial_y^{j-a}b).
\label{eq:toruspersistence}
\end{equation}
The same inverse absorbs the highest derivative. Induction, the trace embedding, and \eqref{eq:toruscoordinate} prove persistence on the original interval. High norms enter the right-hand side, not the lifespan.
Thus \eqref{eq:toruslocaldifference} and persistence hold on the same fixed interval, independently of higher initial norms. This completes the proof.
\end{proof}

\subsection{Sources, short windows, and the uniqueness class}\label{subsec:sources-short-windows-and-the-uniqueness-class}

We define the periodic differentiated source and verify its locality before identifying the physical equation. A zero-trace short-window estimate then extends uniqueness beyond the original contraction ball.
The periodic endpoint product is defined through the source estimate, not as an ordinary square in $X_{\mathbb T}$. From \eqref{eq:torus-local-estimates},
\begin{equation}
\begin{aligned}
\|Q_\lambda(b,b)-Q_\lambda(c,c)\|_{\mathcal F_\lambda}
&\le C_p\lambda^{1/2}(\|b\|_W+\|c\|_W)\|b-c\|_W,\\
|\langle F,\varphi\rangle|
&\lesssim\|F\|_{\mathcal F_\lambda}
\|\langle q\rangle^{1/2}\langle\sigma+q^2\rangle\widetilde\varphi\|_2.
\end{aligned}
\label{eq:torussource}
\end{equation}
In particular $\mathcal F_\lambda\hookrightarrow H_s^{-1}H_y^{-5/2}\hookrightarrow\mathcal D'$.

For rough data the differentiated source must be constructed before the two-field equation is recovered. The next lemma proves that the local source is independent of extension and compatible with smooth approximation.

\begin{lemma}\label{lem:periodic-source}
The source $Q_\lambda$ has a restriction independent of the full-time extension. It agrees with the classical source for smooth functions. The physical differentiated source is
\begin{equation}
\mathcal N(U)(t,x)
=\lambda^5[\partial_y\Omega_\lambda^{-1}Q_\lambda(b,b)]
(\lambda^2t,\lambda x),
\label{eq:torusphysicalsource}
\end{equation}
where $\Omega_\lambda^{-1}$ is zero on the constant mode. It is locally in $H_t^{-1}H_x^{-7/2}$ and is the limit of the smooth differentiated sources.
The recovered zero-mean fields satisfy
\[
u_t=\partial_x v,\qquad v_t=\partial_x (u-\partial_x^2u)-\mathcal N(U)
\]
in distributions.
\end{lemma}

\begin{proof}
We prove locality and agreement with the classical source, establish the regularity and convergence in \eqref{eq:torusphysicalsource}, and then recover the two-field equation using the zero-mean condition.

Finite spatial Fourier truncations and compact time-frequency truncations are dense in $W$ and $\mathcal F_\lambda$. For the dyadic term, truncate its $\ell^1$ tail first. Smooth frequency approximation gives smooth functions rapidly decreasing in time. Time convolution with a compactly supported smooth approximate identity, together with periodic spatial convolution, is uniformly bounded and converges strongly in both spaces: its Fourier multiplier is bounded and converges pointwise to one.

Let $b,c$ agree on an open time interval $I$. Their convolutions agree on every $I'\Subset I$ once the time support is sufficiently small. For smooth convolutions the source is the classical product; this follows by approximation both in $W$ and in a sufficiently high Sobolev space. Its restriction on $I'$ is therefore the same. Passing to the limit by \eqref{eq:torussource} proves locality. The same approximation proves locality of the zero-data Duhamel extension on intervals containing the initial time. No boundedness of an ordinary time cutoff on $W$ is used.

The source has zero spatial mode. At $q=n/\lambda$, $n\ne0$,
\[
\Omega_\lambda(q)^{-1}=1+n^{-2}\le2.
\]
Thus \eqref{eq:torusphysicalsource} is a continuous distributional pullback followed by a bounded nonzero-mode inverse and one derivative. For smooth $U$ it equals $\partial_x (u^2)$; \eqref{eq:torussource} proves convergence for the local fixed-point approximations.

Writing $w=u+i\beta$, where $\beta=(1-\partial_x^2)^{-1}\partial_x v$, gives
\[
u_t=(1-\partial_x^2)\beta=\partial_x v,\qquad
\beta_t=\partial_x^2u+Q_{\rm phys},\qquad
Q_{\rm phys}(t,x)=\lambda^4Q_\lambda(b,b)(\lambda^2t,\lambda x).
\]
Hence
\[
\partial_x \{v_t-\partial_x (u-\partial_x^2u)+\mathcal N(U)\}=0.
\]
Every term inside braces has zero spatial mean. The kernel of $\partial_x$ consists of spatial constants, even for time-distribution coefficients, so the expression in braces vanishes. This proves the two-field equation and excludes an undetermined velocity constant.
The locality, source convergence, and two-field equation prove the lemma.
\end{proof}

Uniqueness in the whole specified restriction class requires more than uniqueness in the contraction ball. We first show that a zero-trace remainder becomes small on a short window, without a general time-cutoff estimate.

\begin{lemma}\label{lem:zero-trace}
If $B\in W$ and $B(0)=0$, then
\begin{equation}
\|B\|_{W((-a,a))}\longrightarrow0\qquad(a\downarrow0).
\label{eq:torusshortwindow}
\end{equation}
\end{lemma}

\begin{proof}
We prove \eqref{eq:torusshortwindow} by approximating the zero-trace function with a smooth zero-trace function and estimating the latter on a shrinking interval.

Choose smooth $V$ with $\|B-V\|_W<\epsilon$. The trace estimate gives
$\|V(0)\|_{H^{-1/2}}\le C\epsilon$. Set
\[
F=V-\psi(s)e^{is\partial_y^2}V(0),\qquad F(0)=0.
\]
For $\chi\in C_c^\infty((-2,2))$, equal to one on $[-1,1]$, let
$\chi_a(s)=\chi(s/a)$ and $g(s)=e^{-is\partial_y^2}F(s)$. Smoothness and the zero trace give
\[
\|\chi_a g\|_{L_s^2H_y^{-1/4}}\le C_Va^{3/2},\qquad
\|\chi_a g\|_{H_s^1H_y^{-1/4}}\le C_Va^{1/2}.
\]
Interpolation and $X_{\rm Sch}^{-1/4,3/4}\hookrightarrow W$ imply
\[
\|B\|_{W((-a,a))}
\le C\epsilon+\|\chi_aF\|_W
\le C\epsilon+C_Va^{3/4}.
\]
First let $a\downarrow0$, then $\epsilon\downarrow0$. This gives \eqref{eq:torusshortwindow}.
This completes the proof.
\end{proof}

\begin{remark}\label{rem:zero-trace-cutoff}
The cutoff in the proof of Lemma~\ref{lem:zero-trace} is applied only to a fixed smooth zero-trace function. The argument does not require boundedness of ordinary time cutoffs on the full resolution space $W$.
\end{remark}

Define $\mathcal K^{\mathbb T}_\lambda(I)$ as follows. A real zero-mean $U\in C(I;X_{\mathbb T})$ belongs to this class if, near every $t_0\in I$,
\[
b_{\lambda,t_0}(s,y)=\lambda^{-2}w(t_0+\lambda^{-2}s,y/\lambda)
\]
belongs to the restriction space $W$ and solves the Duhamel equation associated with \eqref{eq:torusscaling} and Lemma~\ref{lem:periodic-source}. The same scale $\lambda$ is used at every $t_0$.

Subtracting the common free extension reduces the comparison of two local solutions to the preceding zero-trace statement. The contraction estimate then applies on a sufficiently short common window.

\begin{lemma}\label{lem:periodic-unique}
Let $U,V\in\mathcal K^{\mathbb T}_\lambda(I)$ with $0\in I$ and
\[
U(0)=V(0),\qquad\|U(0)\|_{X_{\mathbb T}}\le r_{\rm loc}.
\]
After choosing $r_{\rm loc}$ as in Lemma~\ref{lem:periodic-local}, there is $a>0$ such that
\[
U=V\quad\text{on }(-a,a)\subset I.
\]
This gives uniqueness on common time intervals along the small global trajectory constructed below.
\end{lemma}

\begin{proof}
We make the zero-trace remainders small on a common window, apply the local difference estimate, and then propagate equality along the small trajectory.

Subtract the free extension of the common trace $f$. Lemma~\ref{lem:zero-trace} gives, on $I_a=(-a,a)$,
\[
\limsup_{a\downarrow0}\|b\|_{W(I_a)},\quad
\limsup_{a\downarrow0}\|c\|_{W(I_a)}\le C_p\|f\|.
\]
Using nearly minimal extensions, the local source property, and then taking the infimum gives
\[
\begin{split}
\|b-c\|_{W(I_a)}\le C_p[\lambda^{-2}
+\lambda^{1/2}(\|b\|_{W(I_a)}+\|c\|_{W(I_a)})]
\|b-c\|_{W(I_a)}.
\end{split}
\]
The coefficient is less than one for sufficiently small $a$, by \eqref{eq:toruscontraction} and smallness of the common trace. Hence $b=c$. The bounded inverse \eqref{eq:toruscoordinate} identifies $U$ and $V$. The equality set is relatively open and closed on any common time interval along the constructed small trajectory, which gives the asserted fixed-class uniqueness.
This completes the proof.
\end{proof}

\begin{remark}\label{rem:fixed-class-uniqueness-scope}
Lemma~\ref{lem:periodic-unique} compares solutions in the same fixed rescaling class. It does not identify classes at different scales or assert uniqueness among arbitrary $C_tX_{\mathbb T}$ distributional solutions.
\end{remark}

\subsection{Global continuation and the endpoint limit}\label{subsec:global-continuation-and-the-endpoint-limit}

We combine the invariant spectral barrier with the fixed local lifespan. Smooth continuation precedes the endpoint limit, and the local difference estimate yields finite-time Lipschitz dependence.
Fix $\lambda_{\mathbb T}=\lambda$ from Lemma~\ref{lem:periodic-local}. Choose $r\le r_{\rm loc}$ as in Lemma~\ref{lem:periodic-spectral}, and use the domain $\mathcal D_{\mathbb T}$ of \eqref{eq:periodic-main-domain}. The coercive bounds \eqref{eq:toruscoercive} show that this domain is open in $X_{\mathbb T,0}$ and contains a neighborhood of zero.

\begin{proof}[Proof of Theorem~\ref{thm:circle}]
We first continue the smooth solutions globally using the spectral barrier. We then prove the endpoint approximation and Lipschitz estimate, identify the equation and uniqueness class, and verify the group and persistence assertions.
For smooth $U_0\in\mathcal D_{\mathbb T}$, Lemma~\ref{lem:periodic-local} gives a smooth solution on the fixed interval. Conservation is used only during its existence in $B_r$. Lemma~\ref{lem:periodic-spectral} gives
\begin{equation}
\|U(t)\|^2\le\frac2{c_2}\mathscr S^\mathbb T(U_0)<r^2/4.
\label{eq:torusbarrier}
\end{equation}
There is no first exit from $B_r$. Every restart has the same lifespan $\tau$ and the same high-regularity persistence. Restarting at intervals $\tau/2$ gives a smooth global solution in both time directions. The overlap is identified by Lemma~\ref{lem:periodic-unique}. No high norm determines the restart time.

For general $U_0\in\mathcal D_{\mathbb T}$, take real symmetric Fourier truncations $U_{0,j}\to U_0$. They preserve the zero means. Continuity of $\mathscr S^\mathbb T$ and the strict inequalities imply $U_{0,j}\in\mathcal D_{\mathbb T}$ for all sufficiently large $j$; spectral monotonicity of truncation is not required. All corresponding smooth global paths obey \eqref{eq:torusbarrier}.

Iterate \eqref{eq:toruslocaldifference} on the already existing smooth global paths, in each time direction separately. For $T<\infty$,
\begin{equation}
\sup_{|t|\le T}\|\Phi_{\mathbb T}^tU_{0,j}
-\Phi_{\mathbb T}^tU_{0,h}\|
\le L_{\rm loc}^{\lceil T/\tau\rceil+1}\|U_{0,j}-U_{0,h}\|.
\label{eq:torusglobaldifference}
\end{equation}
Thus the paths converge in $C([-T,T];X_{\mathbb T})$. Interleaving sequences proves independence; compatibility of these limits gives a path on $\mathbb R$. Passing to the limit proves the same difference estimate for every pair of data in $\mathcal D_{\mathbb T}$.

Each spectral functional is continuous, so
\[
\mathscr S^\mathbb T(\Phi_{\mathbb T}^tU_0)=\mathscr S^\mathbb T(U_0),\qquad
A^\mathbb T_\pm(k;\Phi_{\mathbb T}^tU_0)=A^\mathbb T_\pm(k;U_0).
\]
The limit initially lies in the closed $r/2$ ball; the strict sublevel and \eqref{eq:toruscoercive} restore the strict bound and membership in $\mathcal D_{\mathbb T}$.

On every local window, the full-time extensions converge in $W$ by the contraction estimate. Lemma~\ref{lem:periodic-source} gives convergence of their differentiated sources, the two-field equation, and compatibility on overlapping windows. The limit belongs to $\mathcal K^{\mathbb T}_\lambda$.
Lemma~\ref{lem:periodic-unique} identifies time translations and proves
\[
\Phi_{\mathbb T}^{t+s}=\Phi_{\mathbb T}^t\Phi_{\mathbb T}^s,
\qquad \Phi_{\mathbb T}^0=\mathrm{id},\qquad
(\Phi_{\mathbb T}^t)^{-1}=\Phi_{\mathbb T}^{-t}.
\]
Time continuity and \eqref{eq:torusglobaldifference} imply joint continuity. Finally, the persistence estimate on each of finitely many windows preserves $X^{j-1/2}_{\mathbb T}$ for every integer $j\ge0$. The conservation laws, \eqref{eq:torusglobaldifference}, and the group identities prove part (1). Source convergence and fixed-class uniqueness prove part (2), and the final persistence argument proves part (3). This completes the proof of Theorem~\ref{thm:circle}.
\end{proof}

\begin{remark}\label{rem:periodic-construction}
The periodic construction uses the same third-order spectral functional as the line construction, but its endpoint limit is obtained from local stability. In particular, it requires neither a periodic Green-coordinate inverse nor convergence of the periodic auxiliary Hamiltonian flows.
\end{remark}

\paragraph{Acknowledgments.}
The authors thank Benjamin Harrop-Griffiths for sharing his manuscript, for the exchange concerning the independently developed real-line work, and for suggesting the periodic problem.

\paragraph{Funding information.} M. Shan is supported by the National Natural Science Foundation of China grant 12101629. M. Chen is supported by the Guangdong Basic and Applied Basic Research Foundation grant 2025A1515011925. 

\paragraph{Use of automated tools.}
ChatGPT was used in developing arguments, performing symbolic checks, organizing the proofs, and preparing the manuscript. The authors are responsible for all statements and proofs.

\bigskip
\noindent\begin{minipage}{\textwidth}
\hypersetup{urlcolor=black}
{\addresscaps Department of Mathematics, Jinan University\par
Guangzhou 510632, P.R. China\par}
\noindent\textit{Email address:} \href{mailto:mjchen@jnu.edu.cn}{\texttt{mjchen@jnu.edu.cn}}\par

\medskip
{\addresscaps College of Science, Minzu University of China\par
Beijing 100081, P.R. China\par}
\noindent\textit{Email address:} \href{mailto:smj@muc.edu.cn}{\texttt{smj@muc.edu.cn}}\par
\end{minipage}

\end{document}